\documentstyle{elsart}
\input{amssym.def}
\input{amssym.tex}

\journal{Applied Numerical Mathematics}

\newtheorem{teorema}{Theorem}[section]
\newtheorem{propo}[teorema]{Proposition}
\newtheorem{lema}[teorema]{Lemma}

\newtheorem{ej}[teorema]{Example}
\newtheorem{nota}[teorema]{Remark}

\newcommand{\CC}{\Bbb C}
\newcommand{\PP}{\Bbb P}
\newcommand{\ZZ}{\Bbb Z}
\newcommand{\RR}{\Bbb R}

\newcommand{\DD}{\Bbb D}

\newcommand{\TT}{\Bbb T}

\begin{document}
\begin{frontmatter}

\title{A matrix approach to the computation of quadrature formulas on the unit circle\thanksref{Agrad}}
\author{ Mar\'{\i}a Jos\'e Cantero}
\address{Department of Applied Mathematics,  University of
Zaragoza. Calle Maria de Luna 3, 50018 Zaragoza, Spain.}

\author{Ruym\'an Cruz-Barroso  and }
\author{Pablo Gonz\'alez-Vera \corauthref{cor}}

\corauth[cor]{Corresponding author.} \ead{pglez@ull.es}

\address{Department of Mathematical Analysis. La Laguna
University. 38271 La Laguna. Tenerife. Canary Islands. Spain}

\thanks[Agrad]{The work of the first author was partially supported by Gobierno de Arag\'{o}n-CAI, ``Programa Europa de Ayudas a la Investigaci\'{o}n'' and by a research grant from the Ministry of Education and Science of Spain, project MTM2005-08648-C02-01. The work of the
second and the thrid authors was partially supported by the research
project MTM 2005-08571 of the Spanish Government.}

\date{}

\begin{abstract}

In this paper we consider a general sequence of orthogonal Laurent
polynomials on the unit circle and we first study the equivalences
between recurrences for such families and Szeg\H{o}'s recursion and
the structure of the matrix representation for the multiplication
operator in $\Lambda$ when a general sequence of orthogonal Laurent
polynomials on the unit circle is considered. Secondly, we analyze
the computation of the nodes of the Szeg\H{o} quadrature formulas by
using Hessenberg and five-diagonal matrices. Numerical examples
concerning the family of Rogers-Szeg\H{o} $q$-polynomials are also
analyzed.

\vspace{.5cm}

\end{abstract}
\begin{keyword} Orthogonal Laurent polynomials, Szeg\H{o} polynomials, Recurrence relations, Szeg\H{o} quadrature formulas, Rogers-Szeg\H{o} q-polynomials, Hessenberg matrices, five-diagonal matrices.
\end{keyword}
\end{frontmatter}

\section{Introduction}

As it is known, when dealing with the estimation of the integral
$I_{\sigma}(f)=\int_{a}^{b} f(x)d\sigma(x)$, $\sigma(x)$ being a
positive measure on $[a,b]$ by means of an $n$-point
Gauss-Christoffel quadrature rule, $I_n(f)=\sum_{j=1}^{n} A_j
f(x_j)$ such that $I_{\sigma}(P)=I_n(P)$ for any polynomial of
degree up to $2n-1$, the effective computation of the nodes $\{ x_j
\}_{j=1}^{n}$ and weights $\{ A_j \}_{j=1}^{n}$ in $I_n(f)$ has
become an interesting matter of study both numerical and
theoretical. As shown by Gautschi (see \cite{Ga}, \cite{Ga2} or
\cite{Ga3}) among others, here the basic fact is the three-term
recurrence relation satisfied by the sequence of orthogonal
polynomials for the measure $\sigma$ giving rise to certain
tridiagonal matrices (Jacobi matrices) so that the eigenvalues of
the $n$-th principal submatrix coincide with the nodes $\{ x_j
\}_{j=1}^{n}$ i.e., with the zeros of the $n$-th orthogonal
polynomial. Furthermore, the weights $\{ A_j \}_{j=1}^{n}$ can be
easily expressed in terms of the first component of the normalized
eigenvectors.

In this paper, we shall be concerned with the approximate
calculation of integrals of $2\pi$-periodic functions with respect
to a positive measure $\mu$ on $[-\pi,\pi]$ or more generally with
integrals on the unit circle like $I_{\mu}(f)=\int_{-\pi}^{\pi}
f\left( e^{i\theta} \right) d\mu(\theta)$. Here we will also
propose as an estimation for $I_{\mu}(f)$ an $n$-point quadrature
rule $I_n(f)=\sum_{j=1}^{n}\lambda_j f(z_j)$ with distinct nodes
on the unit circle but now imposing exactness not for algebraic
polynomials but trigonometric polynomials or more generally
Laurent polynomials or functions of the form $L(z)=\sum_{j=p}^{q}
\alpha_jz^j$, $\alpha_j \in \CC$, $p$ and $q$ integers with $p
\leq q$. Now, it should be recalled that Laurent polynomials on
the real line were used by Jones and Thron in the early 1980 in
connection with continued fractions and strong moment problems
(see \cite{JNT} and \cite{JT}) and also implicitly in \cite{JH}.
Their study, not only suffered a rapid development in the last
decades giving rise to a theory of orthogonal Laurent polynomials
on the real line (see e.g. \cite{CC}, \cite{CD}, \cite{He},
\cite{OT}, \cite{On} and \cite{Th}), but it was extended to an
ampler context leading to a general theory of orthogonal rational
functions (see \cite{Bu}).

On the other hand, the rapidly growing interest on problems on the
unit circle, like quadratures, Szeg\H{o} polynomials and the
trigonometric moment problem has suggested to develop a theory of
orthogonal Laurent polynomials on the unit circle introduced by
Thron in \cite{Th}, continued in \cite{Jo}, \cite{JC}, \cite{RO} and
where the recent contributions of Cantero, Moral and Vel\'azquez in
\cite{CM}, \cite{Ca} and \cite{CMV3} has meant an important and
definitive impulse for the spectral analysis of certain problems on
the unit circle. Here, it should be remarked that the theory of
orthogonal Laurent polynomials on the unit circle establishes
features totally different to the theory on the real line because of
the close relation between orthogonal Laurent polynomials and the
orthogonal polynomials on the unit circle (see \cite{RP}).

The purpose of this paper is to study orthogonal Laurent polynomials
as well as the analysis and computation of the nodes and weights of
the so-called Szeg\H o quadrature formulas and it is organized as
follows: In Section 2, sequences of orthogonal Laurent polynomials
on the unit circle are constructed. They satisfy certain recurrence
relations which are proved to be equivalent to the recurrences
satisfied by the family of Szeg\H{o} polynomials, as shown in
Section 3. The multiplication operator in the space of Laurent
polynomials with a general ordering previously fixed is considered
in Section 4. This operator plays on the unit circle the fundamental
role in the five-diagonal representation obtained in \cite{Ca}
analogous to the Jacobi matrices in the real line. Our main result
of this section is to prove that this is the minimal representation
from a different point of view than in \cite{CMV3}. In Section 5 a
matrix approach to Szeg\H{o} quadrature formulas in the more
natural framework of orthogonal Laurent polynomials on the unit
circle is analyzed. Finally, we present some illustrative numerical
examples on computation of the nodes and weights of these quadrature
formulas, by way five-diagonal matrices versus Hessenberg matrices
considering the so-called Rogers-Szeg\H{o} weight function in
Section 6.

\section{Orthogonal Laurent polynomials on the unit circle. Preliminary results}

We start this section with some convention for notations and some
preliminary results. We denote by $\TT :=\{ z \in \CC : |z|=1 \}$
and $\DD :=\{ z : |z|<1 \}$ the unit circle and the open disk on
the complex plane respectively. $\PP=\CC[z]$ is the complex vector
space of polynomials in the variable $z$ with complex
coefficients, $\PP_n:=span \{ 1,z,z^3, \dots, z^n \}$ is the
corresponding vector subspace of polynomials with degree less or
equal than $n$ while $\PP_{-1} =\{0\}$ is the trivial subspace.
$\Lambda:=\CC[z,z^{-1}]$ denotes the complex vector subspace of
Laurent polynomials in the variable $z$ and for $m,n \in \ZZ$, $m
\leq n$, we define the vector subspace $\Lambda_{m,n} = span \{
z^m, z^{m+1},\dots, z^n \}$. Also, for a given function $f(z)$ we
define the ``substar-conjugate'' as $f_*(z)=\overline{f\left(
1/\overline{z} \right)}$ whereas for a polynomial $P(z) \in \PP_n$
its reversed (or reciprocal) as
$P^*(z)=z^nP_*(z)=z^n\overline{P\left( 1/\overline{z} \right)}$.

Throughout the paper, we shall be dealing with a positive Borel measure $\mu$ supported on the unit circle $\TT$,
normalized by the condition $\int_{-\pi}^{\pi} d\mu(\theta)=1$ (i.e, a probability measure). As usual,
the inner product induced by $\mu(\theta)$ is given by
$$\langle f,g \rangle_{\mu}=\int_{-\pi}^{\pi} f\left( e^{i\theta} \right) \overline{g\left( e^{i\theta}
\right)} d\mu(\theta).$$

For our purposes, we start constructing a sequence of subspaces of Laurent polynomials $\{ {\cal L}_n \}_{n=0}^{\infty}$ satisfying
$$dim\left( {\cal L}_n \right) = n+1 \;\;\;,\;\; {\cal L}_n \subset {\cal L}_{n+1} \;\;,\;n=0,1,\ldots.$$
This can be done, by taking a sequence $\{ p(n) \}_{n=0}^{\infty}$
of nonnegative integers such that $p(0)=0$, $0 \leq p(n) \leq n$ and
$s(n)=p(n)-p(n-1) \in \{0,1\}$ for $n=1,2,\ldots$. In the sequel, a
sequence $\{ p(n) \}_{n=0}^{\infty}$ satisfying these requirements
will be said a ``generating sequence''. Then, set
\begin{equation}\label{l}
{\cal L}_n = \Lambda_{-p(n),q(n)}= span \left\{ z^j \;:\;-p(n) \leq
j \leq q(n) \right\} \;,\;q(n):=n-p(n).
\end{equation}
Observe that $\{ q(n) \}_{n=0}^{\infty}$ is also a generating sequence and that $\Lambda=\bigcup_{n=0}^{\infty} {\cal L}_n$ if and only if $\lim_{n\rightarrow \infty} p(n) = \lim_{n\rightarrow \infty} q(n) = \infty$. Moreover,
$${\cal L}_{n+1} = \left\{ \begin{array}{lc}
{\cal L}_n \oplus span \{ z^{q(n+1)} \} &if \;s(n+1)=0 \\
\\
{\cal L}_n \oplus span \{ z^{-p(n+1)} \} &if \;s(n+1)=1 \\
\end{array} \right. .$$
In any case, we will say that $\{ p(n) \}_{n=0}^{\infty}$ has
induced an ``ordering'' in $\Lambda$. Sometimes we will need to
define $p(-1)=0$ and hence $s(0)=0$. Now, by applying the
Gram-Schmidt orthogonalization procedure to ${\cal L}_n$, an
orthogonal basis $\{\psi_0(z), \ldots, \psi_n(z) \}$ can be
obtained. If we repeat the process for each $n=1,2,\ldots$, a
sequence $\{ \psi_n(z) \}_{n=0}^{\infty}$ of Laurent polynomials can
be obtained, satisfying
\begin{equation}\label{ortovarphi}
\begin{array}{lcl}
\psi_n(z) \in {\cal L}_n \backslash {\cal L}_{n-1} \;\;,\;n=1,2,\ldots &, &\psi_0(z) \equiv c \neq 0 \\
\\
\langle \psi_n(z),\psi_m(z) \rangle_{\mu}=\kappa_n \delta_{n,m} \;\;,\;\kappa_n > 0 &, &\delta_{n,m}= \left\{
\begin{array}{ccl}
0 &if &n \neq m \\
1 &if &n=m
\end{array}. \right.
\end{array}
\end{equation}
$\{ \psi_n(z) \}_{n=0}^{\infty}$ will be called a ``sequence of
orthogonal Laurent polynomials for the measure $\mu$ and the
generating sequence $\{ p(n) \}_{n=0}^{\infty}$''. It should be
noted that the orders considered by Thron in \cite{Th} (``balanced''
situation), expanding $\Lambda$ in the ordered basis
$$\Lambda_{0,0} \;\;\;,\;\;\;\Lambda_{-1,0} \;\;\;,\;\;\;\Lambda_{-1,1} \;\;\;,\;\;\;\Lambda_{-2,1} \;\;\;,\;\;\;\Lambda_{-2,2} \;\;\;,\;\;\;\Lambda_{-3,2} \;\;\;,\;\;\;\ldots$$
and
$$\Lambda_{0,0} \;\;\;,\;\;\;\Lambda_{0,1} \;\;\;,\;\;\;\Lambda_{-1,1} \;\;\;,\;\;\;\Lambda_{-1,2} \;\;\;,\;\;\;\Lambda_{-2,2} \;\;\;,\;\;\;\Lambda_{-2,3} \;\;\;,\;\;\;\ldots$$
corresponds to $p(n)=E\left[ \frac{n+1}{2} \right]$ and $p(n)=E\left[ \frac{n}{2} \right]$ respectively, where as usual, $E[x]$ denotes the integer part of $x$ (see \cite{Ca}, \cite{RP} and \cite{Cr} for other properties for these particular orderings). In the sequel we will denote by $\{ \phi_n(z) \}_{n=0}^{\infty}$ the sequence of monic orthogonal Laurent polynomials for the measure $\mu$ and the generating sequence $\{ p(n) \}_{n=0}^{\infty}$, that is, when the leading coefficients are equal to 1 for all $n \geq 0$ (coefficients of $z^{q(n)}$ or $z^{-p(n)}$ when $s(n)=0$ or $s(n)=1$ respectively). Moreover, we will denote by $\{ \chi_n(z) \}_{n=0}^{\infty}$ the sequence of orthonormal Laurent polynomials for the measure $\mu$ and the generating sequence $\{ p(n) \}_{n=0}^{\infty}$, i.e. when $\kappa_n=1$ for all $n \geq 0$ in (\ref{ortovarphi}). This sequence is also uniquely determined by assuming that the leading coefficient in $\chi_n(z)$ is positive for each $n \geq 0$.

On the other hand, when taking $p(n)=0$ for all $n=0,1,\ldots$ then
${\cal L}_n=\Lambda_{0,n}=\PP_n$ so that, the $n$-th monic
orthogonal Laurent polynomial coincides with the $n$-th monic
Szeg\H{o} polynomial (see e.g. \cite{Sz}) which will be denoted by
$\rho_n(z)$ for $n=0,1,\ldots$. This means means that $\rho_0(z)
\equiv 1$ and for each $n \geq 1$, $\rho_n(z) \in \PP_n \backslash
\PP_{n-1}$ is monic and satisfies
\begin{equation}\label{rho}
\begin{array}{l}
\langle \rho_n(z),z^s \rangle_{\mu}=\langle \rho_n^{*}(z),z^t \rangle_{\mu}=0 ,\;s=0,1,\ldots,n-1 \;\;,\;t=1,2,\ldots,n  %\;,\\
\\
\langle \rho_n(z),z^n \rangle_{\mu}=\langle \rho_n^{*}(z),1 \rangle_{\mu} > 0.
\end{array}
\end{equation}
Moreover, we will denote by $\{ \varphi_n(z) \}_{n=0}^{\infty}$ the sequence of orthonormal polynomials on the unit circle for $\mu(\theta)$, i.e., satisfying $\parallel \varphi_n(z) \parallel_{\mu}=\langle \varphi_n(z),\varphi_n(z) \rangle_{\mu}^{1/2}=1$ for all $n \geq 0$. This family is uniquely determined by assuming that the leading coefficient in $\varphi_n(z)$ is positive for each $n \geq 0$ and they are related with the family of monic orthogonal polynomials by $\rho_0(z) \equiv \varphi_0(z) \equiv 1$ and $\rho_n(z)=l_n \varphi_n(z)$ with $l_n=\langle \rho_n(z) , \rho_n(z) \rangle_{\mu}^{1/2}$ for all $n \geq 1$.

Explicit expressions for Szeg\"o polynomials are in general not available and in order to compute them we can make use of the following (Szeg\"o) forward recurrence relations (see e.g. \cite{Sz}):
\begin{equation}\label{le}
\begin{array}{lc}
\rho_{0}(z)=\rho_{0}^{*}(z) \equiv 1 & \\
\rho_{n}(z)=z\rho_{n-1}(z)+\delta_{n}\rho_{n-1}^{*}(z) \; &n \geq 1 \\
\rho_{n}^{*}(z)=\overline{\delta_{n}}z\rho_{n-1}(z)+\rho_{n-1}^{*}(z) \; &n \geq 1 \\
\end{array}
\end{equation}
where $\delta_0=1$ and $\delta_n := \rho_n(0)$ for all $n=1,2,\ldots$ are the so-called {\em Schur parameters} ({\em Szeg\"o}, {\em reflection}, {\em Verblunsky} or {\em Geronimus} parameters, see \cite{BS}) with respect to $\mu(\theta)$. Since the zeros of $\rho_n(z)$ lie in $\DD$, they satisfy $|\delta_n|<1$ for $n \geq 1$. Now, if we introduce the sequence of nonnegative real numbers $\{ \eta_n \}_{n=1}^{\infty}$ by
\begin{equation}\label{eta}
\eta_n=\sqrt{1-|\delta_n|^2} \in (0,1] \;\;\;,\;\;n=1,2,\ldots,
\end{equation}
then, a straightforward computation from (\ref{le}) yields
$$\frac{\langle \rho_n(z),\rho_n(z) \rangle_{\mu}}{\langle \rho_{n-1}(z),\rho_{n-1}(z) \rangle_{\mu}} = \eta_n^2$$
and so, a forward recurrence for the family of orthonormal Szeg\H{o} polynomials is given by:

\begin{equation}\label{lenorma}
\begin{array}{lc}
\varphi_{0}(z)=\varphi_{0}^{*}(z) \equiv 1 & \\
\eta_n \varphi_{n}(z)=z\varphi_{n-1}(z)+\delta_{n}\varphi_{n-1}^{*}(z) \; &n \geq 1 \\
\eta_n \varphi_{n}^{*}(z)=\overline{\delta_{n}}z\varphi_{n-1}(z)+\varphi_{n-1}^{*}(z) \; &n \geq 1. \\
\end{array}
\end{equation}

We conclude this section considering the following results proved in \cite{Ca} and \cite{RO}. The first one establishes the relation between the families of orthogonal Laurent polynomials with respect to the generating sequences $\{ p(n) \}_{n=0}^{\infty}$ and $\{ q(n) \}_{n=0}^{\infty}$ whereas the second one states the relation between the family of orthonormal and monic orthogonal Laurent polynomials for the generating sequence $\{ p(n) \}_{n=0}^{\infty}$ and the family of orthonormal and Szeg\H{o} ordinary polynomials. This last one explains how to construct orthogonal Laurent polynomials on the unit circle from the sequence of Szeg\H{o} polynomials. Here it should be remarked that the situation on the real line, i.e. when dealing with sequences of orthogonal Laurent polynomials with respect to a positive measure supported on the real line, is totally different (for details, see e.g. \cite{CD}).

\begin{propo}\label{conex1}
Let $\{\xi_n (z) \}_{n=0}^{\infty}$ be a sequence of orthogonal Laurent polynomials for the measure $\mu$ and the generating sequence $\{q(n)\}_{n=0}^{\infty}$. Then, $\xi_n(z)=\psi_{n*}(z)$ for all $n \geq 0$, $\{\psi_n (z) \}_{n=0}^{\infty}$ being a sequence of orthogonal Laurent polynomials for the measure $\mu$ and the generating sequence $\{p(n)\}_{n=0}^{\infty}$, where $p(n)=n-q(n)$.
\end{propo}
\begin{flushright}
$\Box$
\end{flushright}

\begin{propo}\label{conex2}
The families $\{ \phi_n(z) \}_{n=0}^{\infty}$ and $\{ \chi_n(z) \}_{n=0}^{\infty}$ are the respective sequences of monic orthogonal and orthonormal Laurent polynomials on the unit circle for a measure $\mu$ and the ordering induced by the generating sequence $\{p(n) \}_{n=0}^{\infty}$, if and only if,
\begin{equation}\label{correspond}
\phi_n(z)=\left\{ \begin{array}{ccl}
\frac{\rho_n(z)}{z^{p(n)}} &if &s(n)=0 \\
\\
\frac{\rho_n^*(z)}{z^{p(n)}} &if &s(n)=1
\end{array} \right. \;\;\;\;\;and\;\;\;\;\;
\chi_n(z)=\left\{ \begin{array}{ccl}
\frac{\varphi_n(z)}{z^{p(n)}} &if &s(n)=0 \\
\\
\frac{\varphi_n^*(z)}{z^{p(n)}} &if &s(n)=1
\end{array} . \right.
\end{equation}
\end{propo}
\begin{flushright}
$\Box$
\end{flushright}

\section{Recurrence relations}

\setcounter{equation}{0}

We start with the next result which establishes a three-term recurrence relation for the monic orthogonal and orthonormal families of Laurent polynomials for the measure $\mu$ and certain (balanced) generating sequences,

\begin{propo}\label{equival1}
Consider the familes $\{ \phi_n(z) \}_{n=0}^{\infty}$ and $\{ \tilde{\phi}_n(z) \}_{n=0}^{\infty}$
%$\left( \{ \chi_n(z) \}_{n=0}^{\infty} \right.$ and $\left. \{ \tilde{\chi}_n(z) \}_{n=0}^{\infty} \right)$
 of monic orthogonal Laurent polynomials for the measure $\mu$ and the generating sequences $p(n)=E \left[ \frac{n+1}{2} \right]$ and $p(n)=E \left[ \frac{n}{2} \right]$ respectively. Set
$$A_{n}= \left\{
\begin{array}{crl}
\delta_{n} &if &\;n \;is \;even \\
\overline{\delta_n} &if &\;n \;is \;odd
\end{array}. \right.$$ Then,
%\begin{enumerate}
\begin{equation}\label{eq1} \begin{array}{c}
\phi_{n}(z)=\left( A_n + \overline{A_{n-1}} z^{(-1)^n} \right) \phi_{n-1}(z) + \eta_{n-1}^2 z^{(-1)^n} \phi_{n-2}(z) \;,\;\;n \geq 2 \\
\\
\phi_{0}(z) \equiv 1 \;\;,\;\;\;\;\phi_{1}(z)=\overline{\delta_1} + \frac{1}{z}.
\end{array}
\end{equation}

%\item \begin{equation}\label{eq2} \begin{array}{c}
%\eta_n \cdot \chi_n(z) = \left( A_n + \overline{A_{n-1}} z^{(-1)^n} \right) \chi_{n-1}(z) + \eta_{n-1} z^{(-1)^n} \chi_{n-2}(z) \;\;,\;\;\;\;n \geq 2 \\
%\\
%\chi_{i}(z) = \frac{\phi_{i}(z)}{\parallel \phi_{i}(z) \parallel_{\mu}} \;\;\;\;\;i=0,1.
%\end{array}
%\end{equation}

\begin{equation}\label{eq3} \begin{array}{c}
\tilde{\phi}_{n}(z)=\left( \overline{A_{n}} + A_{n-1} z^{(-1)^{n+1}} \right) \tilde{\phi}_{n-1}(z) + \eta_{n-1}^2 z^{(-1)^{n+1}} \tilde{\phi}_{n-2}(z) \;,\;\;n \geq 2 \\
\\
\tilde{\phi}_{0}(z) \equiv 1 \;\;,\;\;\;\;\tilde{\phi}_{1}(z)=\delta_1 + z.
\end{array}
\end{equation}

%\item \begin{equation}\label{eq4} \begin{array}{c}
%\eta_n \cdot \tilde{\chi}_{n}(z) = \left( \overline{A_n} + A_{n-1} z^{(-1)^{n+1}} \right) \tilde{\chi}_{n-1}(z) + \eta_{n-1} z^{(-1)^{n+1}} \tilde{\chi}_{n-2}(z) \;\;\;\;,\;\;n \geq 2 \\
%\\
%\tilde{\chi}_{i}(z) = \frac{\tilde{\phi}_{i}(z)}{\parallel \tilde{\phi}_{i}(z) \parallel}_{\mu} \;\;\;\;i=0,1.
%\end{array}
%\end{equation}

%\end{enumerate}
\end{propo}
\begin{flushright}
$\Box$
\end{flushright}

These recurrences were initially proved by Thron in \cite{Th} in the context of continued fractions. An alternative proof is given in \cite{Cr} making use of (\ref{le}) and Proposition \ref{conex2}. We will see now that the recurrences (\ref{le}) and (\ref{eq1}) (the same with (\ref{eq3})) are in fact equivalent.

\begin{teorema}\label{equivalole}
The recurrences (\ref{le}) and (\ref{eq1}) are equivalent.
\end{teorema}

{\em Proof}.- As above said, it remains to show that from the relations (\ref{eq1}) we deduce the Szeg\H{o} recurrence (\ref{le}). Indeed, from Proposition \ref{conex2}, (\ref{eq1}) and by taking ``super-star'' conjugation it follows for $k \geq 1$:
\begin{equation}\label{rewritted1}
\begin{array}{lcl}
\rho_{2k}(z) &= &\left( \delta_{2k} + \delta_{2k-1}z \right)\rho_{2k-1}^*(z) + \eta_{2k-1}^2z^2\rho_{2k-2}(z) \\
\\
&= &\delta_{2k}\rho_{2k-1}^*(z) + z\left[ \delta_{2k-1}\rho_{2k-1}^*(z) + \eta_{2k-1}^2z\rho_{2k-2}(z) \right]
\end{array}
\end{equation}
\begin{equation}\label{rewritted2}
\begin{array}{lcl}
\rho_{2k-1}^*(z) &= &\left( \overline{\delta_{2k-1}}z + \overline{\delta_{2k-2}} \right)\rho_{2k-2}(z) + \eta_{2k-2}^2\rho_{2k-3}^*(z) \\
\\
&= &\overline{\delta_{2k-1}}z\rho_{2k-2}(z) + \left[ \overline{\delta_{2k-2}}\rho_{2k-2}(z) + \eta_{2k-2}^2\rho_{2k-3}^*(z) \right]
\end{array}
\end{equation}
\begin{equation}\label{rewritted3}
\begin{array}{lcl}
\rho_{2k}^*(z) &= &\left( \overline{\delta_{2k}}z + \overline{\delta_{2k-1}} \right)\rho_{2k-1}(z) + \eta_{2k-1}^2\rho_{2k-2}^*(z) \\
\\
&= &\overline{\delta_{2k}}z\rho_{2k-1}(z) + \left[ \overline{\delta_{2k-1}}\rho_{2k-1}(z) + \eta_{2k-1}^2\rho_{2k-2}^*(z) \right]
\end{array}
\end{equation}
\begin{equation}\label{rewritted4}
\begin{array}{lcl}
\rho_{2k-1}(z) &= &\left( \delta_{2k-1} + \delta_{2k-2}z \right)\rho_{2k-2}^*(z) + \eta_{2k-2}^2z^2\rho_{2k-3}(z) \\
\\
&= &\delta_{2k-1}\rho_{2k-2}^*(z) + z\left[ \delta_{2k-2}\rho_{2k-2}^*(z) + \eta_{2k-2}^2z\rho_{2k-3}(z) \right].
\end{array}
\end{equation}
Clearly, the proof will be completed if from (\ref{rewritted1})-(\ref{rewritted4}) we deduce that:
\begin{equation}\label{rewritted5}
\rho_{2k-1}(z)=\delta_{2k-1}\rho_{2k-1}^*(z)+\eta_{2k-1}^2z\rho_{2k-2}(z)
\end{equation}
\begin{equation}\label{rewritted8}
\rho_{2k-2}^*(z)=\overline{\delta_{2k-2}}\rho_{2k-2}(z)+\eta_{2k-2}^2\rho_{2k-3}^*(z)
\end{equation}
\begin{equation}\label{rewritted7}
\rho_{2k-1}^*(z)=\overline{\delta_{2k-1}}\rho_{2k-1}(z)+\eta_{2k-1}^2\rho_{2k-2}^*(z)
\end{equation}
\begin{equation}\label{rewritted6}
\rho_{2k-2}(z)=\delta_{2k-2}\rho_{2k-2}^*(z)+\eta_{2k-2}^2z\rho_{2k-3}(z).
\end{equation}
For this purpose, it will be enough to check that (\ref{rewritted5}) is valid since the proof of (\ref{rewritted8})-(\ref{rewritted6}) follows in a similar way. Thus, set
$$R_{2k-1}(z)=\delta_{2k-1}\rho_{2k-1}^*(z)+\eta_{2k-1}^2z\rho_{2k-2}(z) \in \PP_{2k-1}.$$ Hence, $R_{2k-1}(z)=\sum_{j=0}^{2k-1} \alpha_j \rho_j(z).$
By comparision of monomials $z^{2k-1}$ it follows that $\alpha_{2k-1}=1$. Since
$$\alpha_j=\langle R_{2k-1}(z),\rho_j(z) \rangle_{\mu} = \delta_{2k-1} \langle \rho_{2k-1}^*(z),\rho_j(z) \rangle_{\mu} + \eta_{2k-1}^2 \langle z\rho_{2k-2}(z),\rho_j(z) \rangle_{\mu}$$
it follows from the orthogonality conditions for the Szeg\H{o} polynomials that $\alpha_j=0$ for all $j=1,\ldots,2k-2$. Hence, $R_{2k-1}(z)=\rho_{2k-1}(z)+\alpha_0$. Finally, by comparision of monomials $z^0=1$ it follows that $\alpha_0=0$.
\begin{flushright}
$\Box$
\end{flushright}

Moreover, the following recurrence relations were proved in \cite{Ca}:

\begin{propo}\label{equival2}
The family of orthonormal Laurent polynomials for the measure $\mu$ and the generating sequence $p(n)=E \left[ \frac{n}{2} \right]$ satisties
%\begin{equation}\label{five1} \begin{array}{rcl}
%z \left( \begin{array}{c}
%\tilde{\chi}_{0*}(z) \\
%\tilde{\chi}_{1*}(z)
%\end{array} \right) &= &\left( \begin{array}{c}
%-\delta_1 \\
%\eta_1
%\end{array} \right) \tilde{\chi}_{0*}(z) + \left( \begin{array}{cc}
%-\eta_1\delta_2 & \eta_1 \eta_2  \\
%-\overline{\delta_1}\delta_2 & \overline{\delta_1} \eta_2
%\end{array} \right) \left( \begin{array}{c}
%\tilde{\chi}_{1*}(z) \\
%\tilde{\chi}_{2*}(z) \end{array} \right) \;\;\;, \\
%\\
%z \left( \begin{array}{c}
%\tilde{\chi}_{(2n)*}(z) \\
%\tilde{\chi}_{(2n+1)*}(z)
%\end{array} \right) &= &\left( \begin{array}{cc}
%-\eta_{2n}\delta_{2n+1} &-\overline{\delta_{2n}}\delta_{2n+1} \\
%\eta_{2n}\eta_{2n+1} &\overline{\delta_{2n}}\eta_{2n+1}
%\end{array} \right) \left( \begin{array}{c}
%\tilde{\chi}_{(2n-1)*}(z) \\
%\tilde{\chi}_{(2n)*}(z) \end{array} \right) \;+ \\
%\\
%&&\left( \begin{array}{cc}
%-\eta_{2n+1}\delta_{2n+2} & \eta_{2n+1} \eta_{2n+2}  \\
%-\overline{\delta_{2n+1}}\delta_{2n+2} & \overline{\delta_{2n+1}} \eta_{2n+2}
%\end{array} \right) \left( \begin{array}{c}
%\tilde{\chi}_{(2n+1)*}(z) \\
%\tilde{\chi}_{(2n+2)*}(z) \end{array} \right) \;\;,\;\;n \geq 1
%\end{array}
%\end{equation}
%and
\begin{equation}\label{five2} \begin{array}{rcl}
z \tilde{\chi}_0(z) &= &-\delta_1 \tilde{\chi}_0(z) + \eta_1 \tilde{\chi}_1(z) \;\;\;, \\
\\
z \left( \begin{array}{c}
\tilde{\chi}_{2n-1}(z) \\
\tilde{\chi}_{2n}(z)
\end{array} \right) &= &\left( \begin{array}{cc}
-\eta_{2n-1}\delta_{2n} &-\overline{\delta_{2n-1}}\delta_{2n} \\
\eta_{2n-1}\eta_{2n} &\overline{\delta_{2n-1}}\eta_{2n}
\end{array} \right) \left( \begin{array}{c}
\tilde{\chi}_{2n-2}(z) \\
\tilde{\chi}_{2n-1}(z) \end{array} \right) \;+ \\
\\
&&\left( \begin{array}{cc}
-\eta_{2n}\delta_{2n+1} & \eta_{2n} \eta_{2n+1}  \\
-\overline{\delta_{2n}}\delta_{2n+1} & \overline{\delta_{2n}} \eta_{2n+1}
\end{array} \right) \left( \begin{array}{c}
\tilde{\chi}_{2n}(z) \\
\tilde{\chi}_{2n+1}(z) \end{array} \right) \;\;,\;\;n \geq 1
\end{array}
\end{equation}
%respectively.
\begin{flushright}
$\Box$
\end{flushright}
\end{propo}

A similar matrix recurrence is deduced if the generating sequence $p(n)=E \left[ \frac{n+1}{2} \right]$ is considered (see \cite{Ca}). Once we have proved that the Szeg\H{o} recurrence (\ref{le}) is equivalent with both given in Proposition \ref{equival1} we will see now the equivalence with the recurrence given in Proposition \ref{equival2}. It will be proved in the case $p(n)=E\left[ \frac{n}{2} \right]$. The proof when the generating sequence $p(n)=E\left[ \frac{n+1}{2} \right]$ is considered can be established from Proposition \ref{conex1} just by taking ``super-star'' conjugation.

\begin{teorema}\label{iff}
The recurrence relations given in Propositions \ref{equival1} and \ref{equival2} are equivalent.
\end{teorema}

{\em Proof}.- First of all, when the family of orthonormal Laurent polynomials is considered, (\ref{eq3}) becomes
\begin{equation}\label{eq4} \begin{array}{c}
\eta_n \cdot \tilde{\chi}_{n}(z) = \left( \overline{A_n} + A_{n-1} z^{(-1)^{n+1}} \right) \tilde{\chi}_{n-1}(z) + \eta_{n-1} z^{(-1)^{n+1}} \tilde{\chi}_{n-2}(z) \;,\;\;n \geq 2 \\
\\
\tilde{\chi}_{i}(z) = \frac{\tilde{\phi}_{i}(z)}{\parallel \tilde{\phi}_{i}(z) \parallel}_{\mu} \;\;\;\;i=0,1.
\end{array}
\end{equation}
Since the recurrence (\ref{eq4}) is equivalent to (\ref{le}) and the recurrences given in Proposition \ref{equival2} were obtained from this, it remains to prove that (\ref{five2}) implies (\ref{eq4}). From Proposition \ref{conex2} and since the generating sequence $p(n)=E\left[ \frac{n}{2} \right]$ is considered, it follows that (\ref{five2}) is equivalent to the relations
\begin{equation}\label{qq1}
\begin{array}{rcl}
z^2 \varphi_{2n-1}(z) &=
&\eta_{2n}\eta_{2n+1}\varphi_{2n+1}(z)-\eta_{2n}\delta_{2n+1}\varphi_{2n}^*(z)
\\
&&-\overline{\delta_{2n-1}}\delta_{2n}z\varphi_{2n-1}(z)-\eta_{2n-1}\delta_{2n}z\varphi_{2n-2}^*(z)
\end{array}
\end{equation}
\begin{equation}\label{qq2}
\begin{array}{rcl}
z\varphi_{2n}^*(z) &=
&\overline{\delta_{2n}}\eta_{2n+1}\varphi_{2n+1}(z)-\overline{\delta_{2n}}\delta_{2n+1}\varphi_{2n}^*(z)
\\
&&+\overline{\delta_{2n-1}}\eta_{2n}z\varphi_{2n-1}(z)+\eta_{2n-1}\eta_{2n}z\varphi_{2n-2}^*(z).
\end{array}
\end{equation}
Again, from Proposition \ref{conex2} we have to prove from (\ref{qq1}) and (\ref{qq2}) that
\begin{equation}\label{qq3}
\eta_{2n}\varphi_{2n}^*(z) = \left( \overline{\delta_{2n}}z + \overline{\delta_{2n-1}} \right) \varphi_{2n-1}(z) + \eta_{2n-1}\varphi_{2n-2}^*(z)
\end{equation}
and
\begin{equation}\label{qq4}
\eta_{2n+1}\varphi_{2n+1}(z)= \left( \delta_{2n+1} + \delta_{2n}z \right) \varphi_{2n}^*(z) + \eta_{2n}z^2 \varphi_{2n-1}(z)
\end{equation}
holds. By one hand, since $\eta_{2n} \neq 0$ it follows from (\ref{qq1})
$$\begin{array}{rcl}
\eta_{2n}z^2\varphi_{2n-1}(z) &= &\left(1 - |\delta_{2n}|^2 \right) \eta_{2n+1} \varphi_{2n+1}(z) - \left(1 - |\delta_{2n}|^2 \right) \delta_{2n+1}\varphi_{2n}^*(z) \\
\\
&&-\eta_{2n}\delta_{2n}\overline{\delta_{2n-1}}z\varphi_{2n-1}(z)-\eta_{2n}\eta_{2n-1}\delta_{2n}z\varphi_{2n-2}^*(z) \\
\\
&= &\eta_{2n+1}\varphi_{2n+1}(z) - \delta_{2n+1}\varphi_{2n}^*(z) - \delta_{2n} [ \overline{\delta_{2n}}\eta_{2n+1}\varphi_{2n+1}(z) \\
\\
&&- \overline{\delta_{2n}}\delta_{2n+1}\varphi_{2n}^*(z) + \eta_{2n}\overline{\delta_{2n-1}}z\varphi_{2n-1}(z) + \eta_{2n}\eta_{2n-1}z\varphi_{2n-2}^*(z) ]
\end{array}$$
and now from (\ref{qq2}) the relation (\ref{qq4}) is deduced. By other hand, from (\ref{qq2}) and since $\eta_{2n} \neq 0$ it follows that
$$\begin{array}{rcl}
\eta_{2n}z\varphi_{2n}^*(z) &= &\overline{\delta_{2n}}\eta_{2n}\eta_{2n+1}\varphi_{2n+1}(z)-\overline{\delta_{2n}}\delta_{2n+1}\eta_{2n}\varphi_{2n}^*(z) + \overline{\delta_{2n-1}}\eta_{2n}^2z\varphi_{2n-1}(z) \\
\\
&&+\eta_{2n}^2\eta_{2n-1}z\varphi_{2n-2}^*(z) \\
\\
&=& \overline{\delta_{2n}}\eta_{2n} \left[ \eta_{2n+1}\varphi_{2n+1}-\delta_{2n+1}\varphi_{2n}^*(z) \right] + z \times \\
\\
&&\left[ \overline{\delta_{2n-1}}\varphi_{2n-1} -
\overline{\delta_{2n-1}}|\delta_{2n}|^2\varphi_{2n-1} +
\eta_{2n-1}\varphi_{2n-2}^*(z)-\eta_{2n-1}|\delta_{2n}|^2\varphi_{2n-2}^*(z)
\right].
\end{array}$$
Now, from (\ref{qq4})
$$\begin{array}{rcl}
\eta_{2n}\varphi_{2n}^*(z) &= &\overline{\delta_{2n}}\eta_{2n} \left[ \eta_{2n}z \varphi_{2n-1}(z) + \delta_{2n}\varphi_{2n}^*(z) \right] + \overline{\delta_{2n-1}}\varphi_{2n-1} - \overline{\delta_{2n-1}}|\delta_{2n}|^2\varphi_{2n-1}\\
\\
&&+ \eta_{2n-1}\varphi_{2n-2}^*(z)-\eta_{2n-1}|\delta_{2n}|^2\varphi_{2n-2}^*(z) \\
\end{array}$$
which is equivalent to
$$\eta_{2n}^2\left[ \eta_{2n} \varphi_{2n}^*(z)-\overline{\delta_{2n}}z\varphi_{2n-1}(z)-\overline{\delta_{2n-1}}\varphi_{2n-1}(z)-\eta_{2n-1}\varphi_{2n-2}^*(z) \right] = 0.$$
So, (\ref{qq3}) is proved since $\eta_{2n} \neq 0$.

\begin{flushright}
$\Box$
\end{flushright}

\begin{nota}
Observe that the three-term recurrences given in Proposition \ref{equival1} involves multiplication by $z$
and $z^{-1}$ together whereas the recurrences given in Poposition \ref{equival2} involves only multiplication by
$z$. These last relations will play a fundamental role in the next section.
\end{nota}

Untill the present moment we have deduced recurrences for the families of orthogonal Laurent polynomials when
the generating sequences associated with the balanced orderings are considered. In the rest of this section we
will consider an arbitrary generating sequence, starting with a three-term recurrence relation for
$\{ \phi_n(z) \}_{n=0}^{\infty}$ involving multiplication by $z$ and $z^{-1}$. Here, we define $p(-1)=0$ and so $s(0)=0$.

\begin{teorema}\label{ttgeneral}
The family of monic orthogonal Laurent polynomials $\{ \phi_n(z) \}_{n=0}^{\infty}$ with respect to the measure
$\mu$ and the generating sequence $\{ p(n) \}_{n=0}^{\infty}$ satisfies for $n \geq 2$ the three-term
recurrence relation
\begin{equation}\label{tt}
\begin{array}{rcl}
\phi_n(z) &= &\left( A_n B_n + C_n z^{1-2s(n)} \right) \phi_{n-1}(z)
+ \\
&&(-1)^{1 + s(n-2) - s(n-1) } D_n E_n \eta_{n-1}^2 z^{1- s(n)-
s(n-2) } \phi_{n-2}(z)
\end{array}
\end{equation}
with initial conditions
\begin{equation}\label{ttinit}
\phi_0(z) \equiv 1 \;\;,\;\;\;\phi_1(z) = K_1 + z^{1-2s(1)}
%%%%%phi_2(z)=P_2z^{2-s(1)-s(2)}+Q_2z^{1-s(1)-s(2)}+R_2z^{-s(1)-s(2)}%%%%%
\end{equation}
where
%\begin{equation}\label{KP}
$K_1=\left\{ \begin{array}{ccc}
\delta_1 &if &s(1)=0 \\
\overline{\delta_1} &if &s(1)=1
\end{array} \right.$
%\;\;\;,\;\;\;
%P_2=\left\{ \begin{array}{ccc}
%1 &if &s(2)=0 \\
%\overline{\delta_2} &if &s(2)=1
%\end{array} \right.
%\end{equation}
%\begin{equation}\label{QR}
%Q_2=\left\{ \begin{array}{ccc}
%\delta_1+\overline{\delta_1}\delta_2 &if &s(2)=0 \\
%\overline{\delta_1} + \delta_1\overline{\delta_2} &if &s(2)=1
%\end{array} \right. \;\;\;,\;\;\;
%R_2=\left\{ \begin{array}{ccc}
%\delta_2 &if &s(2)=0 \\
%1 &if &s(2)=1
%\end{array} \right.
%\end{equation}
and for $n \geq 2$,
\begin{equation}\label{abn}
A_n= \left\{ \begin{array}{lcl}
1 &if &s(n) \neq s(n-1) \\
\delta_{n-1}^{2s(n)-1} &if &s(n) = s(n-1) \;\;,\;s(n-2)=0 \\
\left( \overline{\delta_{n-1}} \right)^{1-2s(n)} &if &s(n) = s(n-1)
\;\;,\;s(n-2)=1 \end{array} \right. ,
\end{equation}
\begin{equation}\label{bn}
 B_n= \left\{ \begin{array}{lcl}
\delta_n &if &s(n)=0 \\
\overline{\delta_n} &if &s(n)=1
\end{array} \right. ,
\end{equation}
\begin{equation}\label{cn}
C_n= \left\{ \begin{array}{lcl}
1 &if &s(n)=s(n-1) \\
\delta_{n-1}^{s(n-1)-s(n)} &if &s(n) \neq s(n-1) \;\;,\;s(n-2)=0 \\
\left( \overline{\delta_{n-1}} \right)^{s(n)-s(n-1)} &if &s(n) \neq s(n-1) \;\;,\;s(n-2)=1 \\
\end{array} \right. ,
\end{equation}
\begin{equation}\label{dn}
D_n= \left\{ \begin{array}{lcl}
\delta_n &if &s(n-1) = s(n)=0 \\
\overline{\delta_{n}} &if &s(n-1) = s(n) =1 \\
1 &if &s(n-1) \neq s(n) \\
\end{array} \right. ,
\end{equation}
\begin{equation}\label{den}
E_n= \left\{ \begin{array}{lcl}
1/\delta_{n-1} &if &s(n-2) = s(n-1)=0 \\
1/\overline{\delta_{n-1}} &if &s(n-2) = s(n-1)=1 \\
1 &if &s(n-2) \neq s(n-1) \\
\end{array} \right. .
\end{equation}
In the cases $s(n-2)=s(n-1)$ the three-term recurrence relation
holds if and only if $\delta_{n-1} \neq 0$.
\end{teorema}

{\em Proof}.- The initial conditions (\ref{ttinit}) follows from Proposition \ref{conex2} and (\ref{le}). Consider, for $n \geq 2$, the eight cases $\left( s(n),s(n-1),s(n-2) \right) = (i,j,k)$ with $i,j,k \in \{ 0,1 \}$. As we said, the proofs in the balanced situations $(0,1,0)$ and $(1,0,1)$ were given in \cite{Cr} making use of Proposition \ref{conex2} and the Szeg\H{o} recurrence (\ref{le}). Proceeding in the same way in the remaining cases we obtain:
\begin{enumerate}
%\item Case $(0,1,0)$:
%$$\phi_n(z)=\left( \delta_n + \delta_{n-1}z \right) \phi_{n-1}(z) + z\eta_{n-1}^2 \phi_{n-2}(z).$$
%\item Case $(1,0,1)$:
%$$\phi_n(z)=\left( \overline{\delta_n} + \frac{\overline{\delta_{n-1}}}{z} \right) \phi_{n-1}(z) + \frac{1}{z}\eta_{n-1}^2 \phi_{n-2}(z).$$
\item Case $(0,0,0)$:
\begin{equation}\label{000}
\phi_n(z)=\left( \frac{\delta_n}{\delta_{n-1}} + z \right) \phi_{n-1}(z) - \frac{\delta_n}{\delta_{n-1}}\eta_{n-1}^2z \phi_{n-2}(z).
\end{equation}
\item Case $(0,0,1)$:
\begin{equation}\label{001}
\phi_n(z)=\left( \delta_n\overline{\delta_{n-1}} + z \right) \phi_{n-1}(z) + \delta_n\eta_{n-1}^2 \phi_{n-2}(z).
\end{equation}
\item Case $(0,1,1)$:
\begin{equation}\label{011}
\phi_n(z)= \left( \delta_n + \frac{z}{\overline{\delta_{n-1}}} \right) \phi_{n-1}(z) - \frac{1}{\overline{\delta_{n-1}}}\eta_{n-1}^2 \phi_{n-2}(z).
\end{equation}
\item Case $(1,0,0)$:
\begin{equation}\label{100}
\phi_n(z)=\left( \overline{\delta_n} + \frac{1}{\delta_{n-1}z} \right) \phi_{n-1}(z) - \frac{1}{\delta_{n-1}}\eta_{n-1}^2 \phi_{n-2}(z).
\end{equation}
\item Case $(1,1,0)$:
\begin{equation}\label{110}
\phi_n(z)=\left( \overline{\delta_n}\delta_{n-1} + \frac{1}{z} \right) \phi_{n-1}(z) + \overline{\delta_n}\eta_{n-1}^2 \phi_{n-2}(z).
\end{equation}
\item Case $(1,1,1)$:
\begin{equation}\label{111}
\phi_n(z)=\left( \frac{\overline{\delta_n}}{\overline{\delta_{n-1}}} + \frac{1}{z} \right) \phi_{n-1}(z) - \frac{\overline{\delta_n}}{\overline{\delta_{n-1}}z}\eta_{n-1}^2 \phi_{n-2}(z).
\end{equation}
\end{enumerate}
Finally, all these three-term recurrences can be expressed as in (\ref{tt}) along with (\ref{abn})-(\ref{den}) in a whole pattern.
\begin{flushright}
$\Box$
\end{flushright}

\begin{nota}
When the family of orthonormal Laurent polynomials $\{\chi_n(z) \}_{n=0}^{\infty}$ is considered, the relation (\ref{tt}) becomes
\begin{equation}\label{ttnormapn}
\begin{array}{rcl}
\eta_n \chi_n(z) &= &\left( A_n B_n + C_n z^{1-2s(n)} \right)
\chi_{n-1}(z) + \\
&&(-1)^{1 + s(n-2) - s(n-1) } D_n E_n \eta_{n-1} z^{1- s(n)- s(n-2)
} \chi_{n-2}(z)
\end{array}
\end{equation}
with initial conditions $\chi_{i}(z) = \frac{\phi_{i}(z)}{\parallel \phi_{i}(z) \parallel_{\mu}}$ for $i=0,1.$
\end{nota}

In order to complete the equivalences between the recurrence relations it remains to show that the recurrence (\ref{tt})-(\ref{den}) implies the Szeg\H{o} recurrence (\ref{le}). The proof in the balanced situations is given in Theorem \ref{equivalole} and the proof when a general generating sequence is considered can be done with similar arguments starting from the corresponding recurrence (\ref{000})-(\ref{111}). In the next result we consider a general generating sequence and prove a similar result without using orthogonality conditions so that Proposition \ref{conex2} can not be used.

\begin{teorema}\label{sinorto}
Let $\{ p(n) \}_{n \geq 0}$ be a generating sequence and define $p(-1)=0$. Consider $\{ s(n) \}_{n \geq 0}$ defined by $s(n)=p(n)-p(n-1)$ and an arbitrary given sequence of complex numbers $\{ \delta_n \}_{n \geq 0}$ with $\delta_0=1$ and $|\delta_n| \neq 1$ for all $n \geq 1$. Suppose that $\delta_{n-1} \neq 0$ if $s(n-2)=s(n-1)$. Let $\{ \phi_{n}(z) \}_{n \geq 0}$ be the sequence of Laurent polynomials defined by the recurrence relation (\ref{tt}) along with (\ref{abn})-(\ref{den}) and with initials conditions (\ref{ttinit}).
Then,
\begin{enumerate}
\item For all $n \geq 0$, $\phi_n(z)$ is a monic Laurent
polynomial with $\phi_n(z) \in {\cal L}_n \backslash {\cal
L}_{n-1}$. \item Write $\phi_n(z)=\frac{N_n(z)}{z^{p(n)}}$ with
$N_n(z) \in \PP_n$ and set $F_0(z) \equiv N_0(z) \equiv 1$ and
\begin{equation}\label{Fn}
F_n(z)=\left\{ \begin{array}{ccl}
N_n(z) &if &s(n)=0 \\
N_n^*(z) &if &s(n)=1
\end{array} \right.  \;\;\;,\;\;n \geq 1.
\end{equation}
Then, $\{ F_n(z) \}_{n=1}^{\infty}$ satisfies the recurrence relation
\begin{equation}\label{Nn}
\begin{array}{c}
F_{n}(z)=zF_{n-1}(z) +\delta_{n}F_{n-1}^*(z) \\
F_{n}^*(z)=\overline{\delta_n}zF_{n-1}(z) +F_{n-1}^*(z)
\end{array} \;\;\;,\;\;n \geq 1 \;.
\end{equation}
\end{enumerate}
\end{teorema}

{\em Proof}.- $(1)$ Proceeding by induction it is proved for all $n
\geq 1$ that $\phi_n(z) \in {\cal L}_n \backslash {\cal L}_{n-1}$
and that the coefficients of monomials $z^{q(n)}$ and $z^{-p(n)}$
are equal to $1$ when $s(n)=0$ and $s(n)=1$ respectively, that the
coefficient of monomial $z^{-p(n)}$ is equal to $\delta_n$ when
$s(n)=0$ and that the coefficient of monomial $z^{q(n)}$ is equal to
$\overline{\delta_n}$ when $s(n)=1$. The claim holds for $n=1$ from
the initial conditions (\ref{ttinit}). Hence, the proof can be
achieved by induction from (\ref{tt}) by splitting into the eight
cases $\left( s(n),s(n-1),s(n-2) \right) = (i,j,k)$, with $i,j,k \in
\{0,1 \}$.

$(2)$ From (\ref{ttinit}) it clearly follows that $N_0(z) \equiv 1$
and
$$N_1(z)=\left\{ \begin{array}{crl}
\delta_1 + z &if &s(1)=0 \\
\overline{\delta_{1}}z + 1 &if &s(1)=1
\end{array} \right. \;,$$
implying that $F_1(z)=zF_0(z) + \delta_1 F_0^*(z)$.
%For $n=2$ it also follows from (\ref{ttinit}) that
%$$N_2(z)=\left\{ \begin{array}{crl}
%\delta_2 + \left( \delta_1 + \delta_2\overline{\delta_1} \right)z + z^2 &if &s(2)=0 \\
%\\
%1 + \left( \overline{\delta_1} + \overline{\delta_2}\delta_1 \right)z + \overline{\delta_2}z^2 &if &s(2)=1
%\end{array} \right.$$
%and hence $F_2(z)=zF_1(z) + \delta_2 F_1^*(z)$.
 Suppose from the hypothesis of induction that
\begin{equation}\label{Fnmenosuno}
\begin{array}{c}
F_{n-1}(z)=zF_{n-2}(z) +\delta_{n-1}F_{n-2}^*(z) \\
F_{n-1}^*(z)=\overline{\delta_{n-1}}zF_{n-2}(z) +F_{n-2}^*(z)
\end{array} .
\end{equation}
From (\ref{tt}) it follows for $n \geq 2$ that
%\begin{equation}\label{recuordi}
$$\begin{array}{ccl}
\frac{N_n(z)}{z^{p(n)}} &= &\left( A_n B_n + C_n z^{1-2s(n)} \right) \frac{N_{n-1}(z)}{z^{p(n)-s(n)}} + \\
&&(-1)^{1+s(n-2)-s(n-1)}D_n E_n \eta_{n-1}^2 z^{1-s(n)-s(n-2)} \frac{N_{n-2}(z)}{z^{p(n)-\left[ s(n)+s(n-1) \right]}}
\end{array}$$
%\end{equation}
implying the following recurrence relation for the family of ordinary polynomials $\{ N_n(z) \}_{n \geq 2}$:
\begin{equation}\label{recuordi}
\begin{array}{ccl}
N_n(z) &= &\left( A_n B_nz^{s(n)} + C_n z^{1-s(n)} \right) N_{n-1}(z) + \\
&&(-1)^{1+s(n-2)-s(n-1)}D_n E_n \eta_{n-1}^2 z^{1+s(n-1)-s(n-2)} N_{n-2}(z).
\end{array}
\end{equation}
We consider now the four cases $\left( s(n),s(n-1),s(n-2) \right) = (0,j,k)$, with $j,k \in \{0,1 \}$. Then, from (\ref{Fnmenosuno}) and relation (\ref{recuordi}) it follows:

\begin{itemize}
\item Case $(0,0,0)$:
$$\begin{array}{ccl}
N_n(z) &= &\left( \frac{\delta_n}{\delta_{n-1}} + z \right)N_{n-1}(z) - \frac{\delta_n}{\delta_{n-1}}\eta_{n-1}^2zN_{n-2}(z) \\
&= &zN_{n-1}(z) + \frac{\delta_{n}}{\delta_{n-1}} \left[ N_{n-1}(z) - \eta_{n-1}^2zN_{n-2}(z) \right] \\
&= &zN_{n-1}(z) + \frac{\delta_{n}}{\delta_{n-1}} \left[ zN_{n-2}(z) +\delta_{n-1}N_{n-2}^*(z) - \eta_{n-1}^2zN_{n-2}(z) \right] \\
&= &zN_{n-1}(z) + \frac{\delta_{n}}{\delta_{n-1}} \left[ |\delta_{n-1}|^2zN_{n-2}(z) + \delta_{n-1}N_{n-2}^*(z)  \right] \\
&= &zN_{n-1}(z) + \delta_{n} \left[ \overline{\delta_{n-1}}zN_{n-2}(z) + N_{n-2}^*(z)  \right] \\
&= &zN_{n-1}(z) + \delta_{n}N_{n-1}^*(z).
\end{array}$$

\item Case $(0,0,1)$:
$$\begin{array}{ccl}
N_n(z) &= &\left( \overline{\delta_{n-1}}\delta_n + z \right)N_{n-1}(z) + \delta_n\eta_{n-1}^2N_{n-2}(z) \\
&= &zN_{n-1}(z) + \delta_{n} \left[ \overline{\delta_{n-1}}N_{n-1}(z) + \eta_{n-1}^2N_{n-2}(z) \right] \\
&= &zN_{n-1}(z) + \delta_{n} \left[ \overline{\delta_{n-1}}zN_{n-2}^*(z) + |\delta_{n-1}|^2N_{n-2}(z) + \eta_{n-1}^2N_{n-2}(z) \right] \\
&= &zN_{n-1}(z) + \delta_{n} \left[ \overline{\delta_{n-1}}zN_{n-2}^*(z) + N_{n-2}(z) \right] \\
&= &zN_{n-1}(z) + \delta_{n} N_{n-1}^*(z).
\end{array}$$

\item Case $(0,1,0)$:
$$\begin{array}{ccl}
N_n(z) &= &\left( \delta_{n}+\delta_{n-1}z \right)N_{n-1}(z) + \eta_{n-1}^2z^2N_{n-2}(z) \\
&= &\delta_n N_{n-1}(z) + z\left[ \delta_{n-1}N_{n-1}(z) + \eta_{n-1}^2zN_{n-2}(z) \right] \\
&= &\delta_n N_{n-1}(z) + z\left[ |\delta_{n-1}|^2zN_{n-2}(z) + \delta_{n-1}N_{n-2}^*(z) + \eta_{n-1}^2zN_{n-2}(z) \right] \\
&= &\delta_n N_{n-1}(z) + z\left[ zN_{n-2}(z) + \delta_{n-1}N_{n-2}^*(z) \right] \\
&= &\delta_n N_{n-1}(z) + zN_{n-1}^*(z).
\end{array}$$

\item Case $(0,1,1)$:
$$\begin{array}{ccl}
N_n(z) &= &\left( \delta_n + \frac{z}{\overline{\delta_{n-1}}} \right)N_{n-1}(z) - \frac{1}{\overline{\delta_{n-1}}}\eta_{n-1}^2zN_{n-2}(z) \\
&= &\delta_nN_{n-1}(z) + z \frac{1}{\overline{\delta_{n-1}}}\left[ N_{n-1}(z) - \eta_{n-1}^2N_{n-2}(z) \right] \\
&= &\delta_nN_{n-1}(z) + z \frac{1}{\overline{\delta_{n-1}}}\left[ \overline{\delta_{n-1}}zN_{n-2}^*(z) + N_{n-2}(z) - \eta_{n-1}^2 N_{n-2} \right] \\
&= &\delta_nN_{n-1}(z) + z \frac{1}{\overline{\delta_{n-1}}}\left[ \overline{\delta_{n-1}}zN_{n-2}^*(z) + |\delta_{n-1}|^2N_{n-2}(z) \right] \\
&= &\delta_nN_{n-1}(z) + z \left[ zN_{n-2}^*(z) + \delta_{n-1}N_{n-2}(z) \right] \\
&= &\delta_nN_{n-1}(z) + z N_{n-1}^*(z).
\end{array}$$
\end{itemize}
Finally, the proof in the four remaining cases follows from the above proofs and from Proposition \ref{conex1} just by taking super-star conjugation.
\begin{flushright}
$\Box$
\end{flushright}

We conclude this section by considering a Favard-type Theorem for the monic sequence of Laurent polynomials $\{ \phi_n(z) \}_{n=0}^{\infty}$ and a general generating sequence $\{ p(n) \}_{n=0}^{\infty}$. A proof in the ordinary polynomial situation (i.e. $p(n)=0$ for all $n$) is given in \cite{Jo} and an alternative simpler approach in \cite{Er} (see also \cite{Ma}). A proof based on the techniques introduced in the Chihara's book \cite{Ch} in the balanced situations $p(n)=E \left[ \frac{n+1}{2} \right]$ and $p(n)=E \left[ \frac{n}{2} \right]$ is given in \cite{RP}. Anyway, the proof presented here based upon the above equivalences between recurrences is very much simpler, and still valid when a general generating sequence is considered. For our purposes, we briefly recall the concept of orthogonality with respect to a Hermitian linear functional. Indeed, let $\{ \mu_n \}_{n=-\infty}^{\infty}$ be a complex sequence satisfying $\mu_n=\overline{\mu_{-n}}$ for all $n=0,1,2,\ldots$ and denote by $\mu$ the linear functional defined on $\Lambda$ by
$$\mu \left( \sum_{j=p}^{q} \alpha_j z^j \right) := \sum_{j=p}^{q} \alpha_j \mu_{-j} \;\;\;,\;\;\;\alpha_j \in \CC \;\;\;\;\;-\infty < p \leq q < +\infty.$$
In terms of $\mu$ we define a bilinear functional $<\cdot ,\cdot >_{\mu}$ on $\Lambda \times \Lambda$ by
\begin{equation}\label{newproductointerior}
<L,M>_{\mu} := \mu \left( L(z)\overline{M(1/\overline{z})} \right) \;\;\;,\;\;\;L,M \in \Lambda.
\end{equation}

Now (see e.g. \cite{Jo}), it will be said that the functional $\mu$ is quasi-definite if and only if the principal submatrices of the infinite Toeplitz moment matrix associated with $\{ \mu_n \}_{-\infty}^{\infty}$ are nonsingular and positive-definite if the determinants of these matrices are positive.
%Denoting by $\Gamma_n$ the $n$-th Toeplitz determinant associated with $\{ \mu_n \}_{n=-\infty}^{\infty}$
%$$\Gamma_{n}=\left| \begin{array}{ccccc}
%\mu_{0} &\mu_{1} &\mu_{2} &\cdots &\mu_{n} \\
%\mu_{-1} &\mu_{0} &\mu_{1} &\cdots &\mu_{n-1} \\
%\mu_{-2} &\mu_{-1} &\mu_{0} &\cdots &\mu_{n-2} \\
%\vdots &\vdots &\vdots &\cdots &\vdots \\
%\mu_{-n} &\mu_{-n+1} &\mu_{-n+2} &\cdots &\mu_{0} \\
%\end{array}\right| \;\;n=0,1,2,\ldots \;\;(\Gamma_{-1}:=1)$$
%we shall call the functional $\mu$ quasi-definite iff $\Gamma_n \neq 0$ and positive-definite iff $\Gamma_n > 0$ for all %$n=0,1,2,\ldots$.
Quasi-definiteness is a necessary and sufficient condition for the
existence of a family of orthogonal Laurent polynomials with respect
to the linear functional (\ref{newproductointerior}) in the sense
that there exists a sequence $\{ R_n(z) \}_{n=0}^{\infty}$ of
Laurent polynomials satisfying $R_n(z) \in {\cal L}_n \backslash
{\cal L}_{n-1}$, $n=1,2,\ldots$ and $<R_n(z), R_m(z)>_{\mu}=k_n
\delta_{n,m} \;, k_n \neq 0$. On the other hand, if the linear
functional $\mu$ is positive-definite then, the associated linear
functional (\ref{newproductointerior}) is an inner product on
$\Lambda \times \Lambda$ and it holds $<R_n(z),R_m(z)>_{\mu}=k_n
\delta_{n,m}$ with $k_n>0$. When $k_n=1$ for all $n=0,1,2,\ldots$
$\{ R_n(z) \}_{n=0}^{\infty}$ will be called ``orthonormal''.

\begin{teorema}[Favard]\label{favard}
Let $\{ p(n) \}_{n \geq 0}$ be a generating sequence and define $p(-1)=0$ and the sequence $\{ s(n) \}_{n \geq 0}$ by $s(n)=p(n)-p(n-1) \in \{ 0,1 \}$. Consider an arbitrary given sequence of complex numbers $\{ \delta_n \}_{n \geq 0}$ with $\delta_0=1$, $|\delta_n| \neq 1$ for all $n \geq 1$ and set $\{ \lambda_n \}_{n \geq 1}$ defined by $\lambda_n=1-|\delta_n|^2$. Suppose also the restriction $\delta_{n-1} \neq 0$ in the cases $s(n-2)=s(n-1)$. Let $\{ \phi_{n}(z) \}_{n \geq 0}$ be the sequence of Laurent polynomials defined for $n \geq 2$ by the recurrence relation
\begin{equation}\label{ttlambda}
\begin{array}{rcl}
\phi_n(z) &= &\left( A_n B_n + C_n z^{1-2s(n)} \right) \phi_{n-1}(z)
\\
&&+ (-1)^{1 + s(n-2) - s(n-1) } D_n E_n \lambda_{n-1} z^{1- s(n)-
s(n-2) } \phi_{n-2}(z)
\end{array}
\end{equation}
with initial conditions (\ref{ttinit}) and where $A_n$, $B_n$, $C_n$, $D_n$ and $E_n$ are complex constants given
by (\ref{abn})-(\ref{den}). Then, for a fixed $\mu_0 \in \RR \backslash \{0\}$
there exists a unique quasi-definite linear functional $\mu$ such that $\mu(1)=\mu_0$
and $\{ \phi_{n}(z) \}_{n=0}^{\infty}$ is the sequence of monic orthogonal Laurent polynomials
with respect to $\mu$ and the ordering induced in $\Lambda$ by the generating sequence $\{ p(n) \}_{n \geq 0}$.
Furthermore, if we take $\mu_0 >0$ then $\mu$ is positive definite if and only if $|\delta_{n}|<1$ for all $n=1,2,\dots$.
\end{teorema}

{\em Proof}.- From Theorem \ref{sinorto} one sees that for $n \geq
1$, $\phi_n(z) \in {\cal L}_n \backslash {\cal L}_{n-1}$ so that
one can write $\phi_n(z)=\frac{N_n(z)}{z^{p(n)}}$ with $N_n(z) \in
\PP_n$. If we define now the sequence of ordinary polynomials $\{
\rho_n(z) \}_{n \geq 0}$ by
$$\rho_n(z)= \left\{ \begin{array}{ccc}
N_n(z) &if &s(n)=0 \\
N_n^*(z) &if &s(n)=1
\end{array} \right. \;,$$
then, the relations (\ref{Nn})-(\ref{Fn}) become (\ref{le}). Thus, making use of the Favard's theorem in the ordinary polynomial situation (see e.g. \cite{Er} or \cite{Jo}) we see that the family $\{ \rho_n(z) \}_{n \geq 0}$ represents the family of monic orthogonal polynomials (Szeg\H{o} polynomials) with respect to a unique quasi-definite linear functional $\mu$ with $\mu(1)=\mu_0$ and positive definite if and only if $\mu_0 >0$ and $|\delta_{n}|<1$ for all $n=1,2,\dots$. Finally, proof follows from Proposition \ref{conex2} that still clearly holds when the measure $\mu$ is replaced for a quasi-definite (positive) linear functional.
\begin{flushright}
$\Box$
\end{flushright}
\section{The multiplication operator in $\Lambda$}

\setcounter{equation}{0}

Through this section it will play a fundamental role the multiplication operator defined on $\Lambda$, namely
$$\begin{array}{rccc}
M: &\Lambda &\longrightarrow &\Lambda \\
&L(z) &\rightarrow &zL(z)
\end{array}.$$
As we have seen, if we consider the sequence of orthonormal
Laurent polynomials with respect to the measure $\mu$ and the
generating sequence $p(n)=0$ for all $n \geq 0$ then the $n$-th
orthonormal Laurent polynomial coincides with the $n$-th
orthonormal Szeg\H{o} polynomial, for all $n \geq 0$. Since the
operator $M$ lets $\PP$ invariant, taking $\{ \varphi_n(z)
\}_{n=0}^{\infty}$ as a basis of $\PP$ then the following matrix
representation of the restriction of $M$ to $\PP$ is obtained (see
e.g. \cite{Go}, \cite{GR} or \cite{BS}):
\begin{equation}\label{hess}
{\cal H}(\delta)= \left( \begin{array}{cccccc}
h_{0,0} &h_{0,1} &0 &0 &0 &\cdots \\
h_{1,0} &h_{1,1} &h_{1,2} &0 &0 &\cdots \\
h_{2,0} &h_{2,1} &h_{2,2} &h_{2,3} &0 &\cdots \\
\vdots &\vdots &\vdots &\vdots &\vdots &\ddots \\
\end{array} \right) .
\end{equation}
The elements of this matrix with Hessenberg structure are given by
\begin{equation}\label{elementhess}
h_{i,j}= \left\{ \begin{array}{ll}
-\overline{\delta_j} \delta_{i+1} \prod_{k=j+1}^{i} \eta_k &if \;\;j=0,1,\ldots,i-1 \;, \\
-\overline{\delta_i} \delta_{i+1} &if \;\;j=i \;, \\
\eta_{i+1} &if \;\;j=i+1 \;.
\end{array} \right.
\end{equation}
If we consider now the generating sequence $p(n)=E \left[ \frac{n}{2} \right]$ it follows from the recurrence given in Proposition \ref{equival2} the following five-diagonal matrix (CMV representation, see \cite{BS}) for the multiplication operator M (see \cite{Ca}) wich can be also expressed as a product of two tri-diagonal ones:
\begin{equation}\label{cmvmatrix}
\begin{array}{ccl}
{\cal C}(\delta) &= &\left( \begin{array}{cccccccc}
-\delta_1 &\eta_1 &0 &0 &0 &0 &0 &\cdots \\
-\eta_1 \delta_2 &-\overline{\delta_1}\delta_2 &-\eta_2 \delta_3 &\eta_2 \eta_3 &0 &0 &0 &\cdots \\
\eta_1 \eta_2 & \overline{\delta_1} \eta_2 &-\overline{\delta_2}\delta_3 &\overline{\delta_2}\eta_3 &0 &0 &0 &\cdots \\
0 &0 &-\eta_3 \delta_4 &-\overline{\delta_3}\delta_4 &-\eta_4 \delta_5 &\eta_4 \eta_5 &0 &\cdots \\
0 &0 &\eta_3 \eta_4 &\overline{\delta_3}\eta_4 &-\overline{\delta_4} \delta_5 &\overline{\delta_4} \eta_5 &0 &\cdots \\
\vdots &\vdots &\vdots &\vdots &\vdots &\vdots &\vdots &\ddots
\end{array} \right) \\
\\
&= &\left( \begin{array}{ccccccc}
1 &0 &0 &0 &0 &0 &\cdots \\
0 &-\delta_2 &\eta_2 &0 &0 &0 &\cdots \\
0 &\eta_2 &\overline{\delta_2} &0 &0 &0 &\cdots \\
0 &0 &0 &-\delta_4 &\eta_4 &0 &\cdots \\
0 &0 &0 &\eta_4 &\overline{\delta_4} &0 &\cdots \\
\vdots &\vdots &\vdots &\vdots &\vdots &\vdots &\ddots
\end{array} \right) \left( \begin{array}{ccccccc}
-\delta_1 &\eta_1 &0 &0 &0 &0 &\cdots \\
\eta_1 &\overline{\delta_1} &0 &0 &0 &0 &\cdots \\
0 &0 &-\delta_3 &\eta_3 &0 &0 &\cdots \\
0 &0 &\eta_3 &\overline{\delta_3} &0 &0 &\cdots \\
0 &0 &0 &0 &-\delta_5 &\eta_5 &\cdots \\
\vdots &\vdots &\vdots &\vdots &\vdots &\vdots &\ddots
\end{array} \right) .
\end{array}
\end{equation}
Moreover, it is easy to check that the matrix representation when the generating sequence $p(n)=E \left[ \frac{n+1}{2} \right]$ is considered is ${\cal C}(\delta)^{T}$.

The aim of this section is to analyze the structure of the matrix representation for the multiplication operator $M$ when an arbitrary generating sequence is considered. We start with the following result which is a generalization of Proposition 2.4 in \cite{Ca} (here the generating sequences associated with the balanced orderings are previously fixed):

\begin{teorema}\label{tk}
Let $\{ \chi_n(z) \}_{n=0}^{\infty}$ be the sequence of orthonormal Laurent polynomials for the measure $\mu$ and the generating sequence $\{ p(n) \}_{n=0}^{\infty}$ and suppose that $\lim_{n \rightarrow \infty} q(n)=\infty$. Then, for each $n \geq 0$ there exists $k=k(n) \geq 1$ and $t=t(n) \geq 1$ such that $z \chi_n(z) \in span \{ \chi_{n-t}(z), \cdots , \chi_{n+k}(z) \}$, i.e.
$$z \chi_n(z) = \sum_{s=n-t}^{n+k} a_{n,s} \chi_s(z) \;\;\;\;,\;\;\;a_{n,s}=\langle z\chi_n(z),\chi_s(z) \rangle_{\mu}.$$
Moreover, $k=k(n)$ and $t=t(n)$ are defined as follows:
\begin{enumerate}
\item $k=1$ if $s(n+1)=0$ and otherwise $k \geq 2$ is defined satisfying $s(n+1)=\cdots=s(n+k-1)=1, s(n+k)=0$.
\item $t=1$ if $s(n-1)=1$ and otherwise $t \geq 2$ is defined satisfying $s(n-1)=\cdots=s(n+1-t)=0, s(n-t)=1$.
\end{enumerate}
\end{teorema}

{\em Proof}.- On the one hand, since $\chi_n(z) \in {\cal L}_n$ then
$$z\chi_n(z) \in z{\cal L}_n = span \{ \frac{1}{z^{p(n)-1}} , \cdots
, z^{q(n)+1} \} \subset {\cal L}_{n+k}$$ with $k=k(n) \geq 1$.
Observe that the existence of $k$ is guaranteed from the condition
$\lim_{n \rightarrow \infty} q(n)=\infty$. On the other hand, since
$\chi_n(z) \perp {\cal L}_{n-1}$ then
$$z\chi_n(z) \perp z{\cal
L}_{n-1} = span \{ \frac{1}{z^{p(n-1)-1}} , \cdots , z^{q(n-1)+1} \}
\supset {\cal L}_{n-1-t}$$ with $t=t(n) \geq 1$. Since $z\chi_n(z)
\subset {\cal L}_{n+k}$ and $z\chi_n(z) \perp {\cal L}_{n-1-t}$ the
proof follows.
\begin{flushright}
$\Box$
\end{flushright}

From Theorem \ref{tk} we analyze now the matrix representation for the operator $M$ with respect to a generating sequence with minimal number of diagonals. As we have seen, the balanced orderings gives rise to five-diagonal matrices and we must ask if there exists other possible generating sequences which gives rise to a matrix representation with less or equal number of diagonals than five. In this respect we start remarking that a five-diagonal representation is obtained when it happens

\begin{equation}\label{5conditions}
\begin{array}{cl}
1. &k(n)=1,\;t(n) \leq 3 \;\;for \;all \;n \\
2. &k(n) \leq 2,\;t(n) \leq 2 \;\;for \;all \;n \\
3. &k(n) \leq 3,\;t(n) = 1 \;\;for \;all \;n. \\
\end{array}
\end{equation}

Hence, we have the following considerations:
\begin{enumerate}
\item
\begin{lema}\label{minilema}
A five-diagonal representation is not obtained if $\{s(n) \}_{n \geq 1}$ contains three or more consecutive zeros or ones.
\end{lema}

{\em Proof}.- From Theorem \ref{tk} it follows that a block $(s(n),s(n+1),s(n+2),s(n+3),s(n+4))=(1,0,0,0,1)$ implies $(k(n),k(n+1),k(n+2))=(1,1,1)$ with $k(n+3) \geq 2$ and $(t(n+1),t(n+2),t(n+3),t(n+4))=(1,2,3,4)$ whereas a block $(s(n),s(n+1),s(n+2),s(n+3),s(n+4))=(0,1,1,1,0)$ implies $(k(n),k(n+1),k(n+2),k(n+3))=(4,3,2,1)$ and $t(n+1) \geq 2$ with $(t(n+2),t(n+3),t(n+4))=(1,1,1)$. The proof follows since the condition (\ref{5conditions}) fulfill in both situations. From this, the proof when $\{s(n) \}_{n \geq 1}$ contains more than three consecutive zeros or ones is trivial.
\begin{flushright}
$\Box$
\end{flushright}
As a consequence of Lemma \ref{minilema}, the number of consecutive zeros or ones in the sequence $\{s(n) \}_{n \geq 1}$ is at most 2.

\item If the number of consecutive zeros or ones is just one for all $n$ then it corresponds to the generating sequences
$$p(n)=E\left[ \frac{n}{2} \right] \;\;\;,\;\;\;p(n)=E\left[ \frac{n+1}{2} \right].$$
and the five-diagonal matrix representation ${\cal C}(\delta)$ and ${\cal C}(\delta)^{T}$ are obtained respectively.
\end{enumerate}

Now, let us concentrate what happens if two consecutive zeros or
ones appears in the sequence $\{ s(n) \}_{n \geq 1}$. Indeed, a
block of the form$$\left( s(n),s(n+1),s(n+2),s(n+3) \right) =
(1,0,0,1)$$ implies $t(n+2)=2$, $k(n+2) \geq 2$, $t(n+3)=3$, $k(n+3)
\geq 1$ and hence a five-diagonal representation is not obtained
since condition (\ref{5conditions}) does not hold. Suppose now a
block of the form $$\left( s(n),s(n+1),s(n+2),s(n+3) \right) =
(0,1,1,0).$$ If $s(n-1)=0$ then $t(n+1)=3$ and $k(n+1)=2$ whereas if
$s(n-1)=1$ then $t(n)=1$, $k(n)=3$ and $t(n+1)=k(n+1)=2$, implying
in both cases that the condition (\ref{5conditions}) is not
satisfied. Observe that this argument is valid for all $n \geq 2$,
but it really holds for $n \geq 0$. Indeed, if we consider the
generating sequences $p(n)=E \left[ \frac{n+1}{2} \right]$ for all
$n \geq 2$ with $p(0)=p(1)=0$ or $p(n)=E \left[ \frac{n}{2} \right]$
for all $n \geq 2$ with $p(0)=0$ and $p(1)=1$ then it is easy to
check that a non-five diagonal matrix representation is obtained.
Summarizing, we can enunciate:

%If $s(n-1)=0$ then $s(n-2)=1$ (since three consecutive zeros are not available) and the block $$\left( s(n-2),s(n-1),s(n),s(n+1) \right) = (1,0,0,1)$$ is obtained, which yields a $l$-diagonal matrix representation with $l \geq 6$. If $s(n-1)=1$ we are considering the block $$\left( s(n-1),s(n),s(n+1),s(n+2),s(n+3) \right) = (1,0,1,1,0)$$
%and hence from Theorem \ref{tk} if holds that $t(n+1)=2, \;k(n+1)=2, \;t(n)=1, \;k(n)=3$, implying again that the condition (\ref{5conditions}) fulfill. Hence, we have proved the following

\begin{teorema}\label{optimal}
The matrix representation for the multiplication operator $M$ is a five-diagonal matrix if and only if $p(n)=E\left[ \frac{n}{2} \right]$ or $p(n)=E\left[ \frac{n+1}{2} \right]$. Moreover, this representation is the narrowest one in the sense that any matrix representation for another different generating sequence gives rise to a $l$-diagonal matrix representation with $l \geq 6$.
\begin{flushright}
$\Box$
\end{flushright}
\end{teorema}

\begin{nota}
The proof of Theorem \ref{optimal} was obtained using orthogonality conditions. This result has been also recently proved in \cite{CMV3} by using operator theory techniques.
\end{nota}

\begin{ej}
Suppose that we expand $\Lambda$ in the ordered basis
$$\{1,z,z^{-1},z^2,z^{-2},z^{-3},z^3,z^4,z^{-4},z^{-5},z^5,z^6,z^{-6},z^{-7},\ldots
\}.$$ Here $\{ s(n) \}_{n \geq 1} = \{
0,1,0,1,1,0,0,1,1,0,0,1,1,\ldots \}$ and $\lim_{n \rightarrow
\infty} q(n)=\infty$. Observe that $z\chi_0(z) \in span\{
\chi_0(z),\chi_1(z) \}$ and $z\chi_1(z) \in span\{
\chi_0(z),\chi_1(z), \chi_2(z), \chi_3(z) \}$. From the sequence $\{
s(n) \}_{n \geq 1}$ and Theorem \ref{tk} we construct the sequences
$\{ t(n) \}_{n \geq 2} = \{ 2,1,2,1,1,2,3,1,1,\ldots \}$ and $\{
k(n) \}_{n \geq 2} = \{1,3,2,1,1,3,2,1,1,\ldots \}$. These sequences
indicate us the number of nonzero elements in each file of the
matrix. Now, the coefficients of the matrix can be obtained in terms
of the orthonormal Szeg\H{o} polynomials by using the formula given
in Theorem \ref{tk} and Proposition \ref{conex2}, resulting
expressions of the form $\langle z^l f(z),g(z) \rangle_{\mu}$ where
$l \geq 0$ and $f(z),g(z) \in \{ \varphi_n(z) \}_{n \geq 0} \cup \{
\varphi_n^*(z) \}_{n \geq 0}$. Some of them are explicitly
calculated in terms of the sequence of Schur parameters (see
\cite[Ch. 1]{BS}). The remainder quantities can be obtained from
these formulas and from the relations
$$\begin{array}{rcl}
\langle z^t \varphi_m^*(z),\varphi_n(z) \rangle_{\mu} &= &\langle z^{t+m-n} \varphi_n^*(z),\varphi_m(z) \rangle_{\mu} \\
\langle z^t \varphi_m^*(z),\varphi_n^*(z) \rangle_{\mu} &= &\langle z^{t+m-n} \varphi_n(z),\varphi_m(z) \rangle_{\mu}
\end{array} \;\;\;,\;\;t \in \ZZ \;\;,\;n,m \geq 0,$$
which holds since $\varphi_n^*(z)=z^n\overline{\varphi_n}(z)$ when $z \in \TT$. The result is the following matrix representation with seven-diagonal structure:
\begin{center}
{\small
$$\left( \begin{array}{cccccccccccc}
-\delta_1 &\eta_1 &0 &0 &0 &0 &0 &0 &0 &0 &0 &\cdots \\
-\eta_1\delta_2 &-\overline{\delta_1}\delta_2 &-\eta_2\delta_3 &\eta_2\eta_3 &0 &0 &0 &0 &0 &0 &0 &\cdots \\
\eta_1\eta_2 &\overline{\delta_1}\eta_2 &-\overline{\delta_2}\delta_3 &\overline{\delta_2}\eta_3 &0 &0 &0 &0 &0 &0 &0 &\cdots \\
0 &0 &-\eta_3\delta_4 &-\overline{\delta_3}\delta_4 &-\eta_4\delta_5 &-\eta_4\eta_5\delta_6 &\eta_4\eta_5\eta_6 &0 &0 &0 &0 &\cdots \\
0 &0 &\eta_3\eta_4 &\overline{\delta_3}\eta_4 &-\overline{\delta_4}\delta_5 &-\overline{\delta_4}\eta_5\delta_6 &\overline{\delta_4}\eta_5\eta_6 &0 &0 &0 &0 &\cdots \\
0 &0 &0 &0 &\eta_5 &-\overline{\delta_5}\delta_6 &\overline{\delta_5}\eta_6 &0 &0 &0 &0 &\cdots \\
0 &0 &0 &0 &0 &-\eta_6\delta_7 &-\overline{\delta_6}\delta_7 &\eta_7 &0 &0 &0 &\cdots \\
0 &0 &0 &0 &0 &-\eta_6\eta_7\delta_8 &-\overline{\delta_6}\eta_7\delta_8 &-\overline{\delta_7}\delta_8 &-\eta_8\delta_9 &-\eta_8\eta_9\delta_{10} &\eta_8\eta_9\eta_{10} &\cdots \\
0 &0 &0 &0 &0 &\eta_6\eta_7\eta_8 &\overline{\delta_6}\eta_7\eta_8 &\overline{\delta_7}\eta_8 &-\overline{\delta_8}\delta_9 &-\overline{\delta_8}\eta_9\delta_{10} &\overline{\delta_8}\eta_9\eta_{10} &\cdots \\
0 &0 &0 &0 &0 &0 &0 &0 &\eta_9 &-\overline{\delta_9}\delta_{10} &\overline{\delta_9}\eta_{10} &\cdots \\
0 &0 &0 &0 &0 &0 &0 &0 &0 &-\eta_{10}\delta_{11} &-\overline{\delta_{10}}\delta_{11} &\cdots \\
\vdots &\vdots &\vdots &\vdots &\vdots &\vdots &\vdots &\vdots &\vdots &\vdots &\vdots &\ddots
\end{array} \right) .$$}
\end{center}
\end{ej}

Finally we conclude this section considering briefly the inverse multiplication operator defined on $\Lambda$, namely
$$\begin{array}{rccc}
N: &\Lambda &\longrightarrow &\Lambda \\
&L(z) &\rightarrow &\frac{1}{z}L(z)
\end{array}.$$
The first result is an analogous to Theorem \ref{tk} which can be proved in the same way:
\begin{teorema}\label{tkbis}
Let $\{ \chi_n(z) \}_{n=0}^{\infty}$ be the sequence of orthonormal Laurent polynomials for the measure $\mu$ and the generating sequence $\{ p(n) \}_{n=0}^{\infty}$ and suppose that $\lim_{n \rightarrow \infty} p(n)=\infty$. Then, for each $n \geq 0$ there exists $\tilde{k}=\tilde{k}(n) \geq 1$ and $\tilde{t}=\tilde{t}(n) \geq 1$ such that $\frac{1}{z} \chi_n(z) \in span \{ \chi_{n-\tilde{t}}(z), \cdots , \chi_{n+\tilde{k}}(z) \}$, i.e.
\begin{equation}\label{tildean}
\frac{1}{z} \chi_n(z) = \sum_{s=n-\tilde{t}}^{n+\tilde{k}} \tilde{a}_{n,s} \chi_s(z) \;\;\;\;,\;\;\;\tilde{a}_{n,s}=\overline{a_{s,n}}.
\end{equation}
Moreover, $\tilde{k}=\tilde{k}(n)$ and $\tilde{t}=\tilde{t}(n)$ are defined as follows:
\begin{enumerate}
\item $\tilde{k}=1$ if $s(n+1)=1$ and otherwise $\tilde{k} \geq 2$ is defined satisfying $s(n+1)=\cdots=s(n+\tilde{k}-1)=0, s(n+\tilde{k})=1$.
\item $\tilde{t}=1$ if $s(n-1)=0$ and otherwise $\tilde{t} \geq 2$ is defined satisfying $s(n-1)=\cdots=s(n+1-\tilde{t})=1, s(n-\tilde{t})=0$.
\end{enumerate}
\end{teorema}
\begin{flushright}
$\Box$
\end{flushright}

Now, proceeding as before, a similar result to Theorem \ref{optimal} can be deduced for the inverse multiplication operator $N$. Furthermore, from (\ref{tildean}) it follows that the matrix representation of the operator $N$ when dealing with the balanced generating sequence $p(n)=E \left[ \frac{n}{2} \right]$ is ${\cal C}(\delta)^*=\overline{{\cal C}(\delta)}^{T}$.

\section{Quadrature formulas on the unit circle}

\setcounter{equation}{0}

Throughout this section we shall be concerned with the estimation of integrals on $\TT$ of the form,
\begin{equation}\label{int}
I_{\mu}(f)= \int_{-\pi}^{\pi} f\left( e^{i\theta} \right) d\mu(\theta).
\end{equation}
As usual, estimations of $I_{\mu}(f)$ may be produced when replacing in (\ref{int}), $f(z)$ by an appropriate
approximating (interpolating) function $L(z)$ so that $I_{\mu}(L)$ can be now easily computed. It seems reasonable to choose as an approximation to $f(z)$ in (\ref{int}) some appropriate Laurent polynomial, because of the density of $\Lambda$ in $C(\TT)=\{ f: \TT \rightarrow \CC \;,\; f \;{\rm continuous} \}$ with respect to the uniform norm (see e.g. \cite{Ra} and \cite{Sz}). So, the so-called \lq\lq quadrature formulas on the unit circle" arise. Indeed, given the integral $I_{\mu}(f)$, by an $n$-point quadrature formula on $\TT$ we mean an expression like
\begin{equation}\label{quad}
I_{n}(f)=\sum_{j=1}^{n} \lambda_j f(z_j) \;\;\;,\;\;z_i \neq z_j \;\;,\;\;i \neq j
\;\;\;,\;\; z_j \in \TT \;\;,\;j=1, \ldots, n,
\end{equation}
where the nodes $\{z_j \}_{j=1}^{n}$ and the coefficients or weights $\{ \lambda_j \}_{j=1}^{n}$ are chosen so that $I_n(f)$ exactly integrates $I_{\mu}(f)$ in subspaces of $\Lambda$ with dimension as large as possible i.e.,
 $I_n(L)=I_{\mu}(L)$ for any $L \in  \Lambda_{-p,q}$ with $p$ and $q$ nonnegative integers depending on $n$ with sum as great as possible. If we first try with subspaces of the form $ \Lambda_{-p,p}$, it can be easily checked that there can not exist an $n$-point quadrature formula $I_n(f)$ as (\ref{quad}) to be exact in $ \Lambda_{-n,n}$. Hence, it holds that  $p \leq n-1$. In \cite{RO} the following  \lq\lq necessary condition" on the nodal polynomial is proved:

\begin{teorema}\label{21}
For $n \geq 1$, let $I_n(f)=\sum_{j=1}^{n} \lambda_j f(z_j)$ with
$z_j \in \TT$, $j=1,\ldots,n$ be exact in $ \Lambda_{-(n-1),n-1}$,
and set $P_n(z)=\prod_{j=1}^{n} (z-z_j)$. Then,
\begin{equation}\label{paranewbis}
P_n(z)=C_n\left[ \rho_n(z) + \tau_n \rho_n^*(z) \right]
\;\;\;,\;\;|\tau_n| =1 \;\;\;,\;\;C_n=\left( 1 + \tau_n
\overline{\delta_n} \right)^{-1}
\end{equation}
\end{teorema}
\begin{flushright}
$\Box$
\end{flushright}
Moreover, in \cite{Jo} it is proved the following converse result (sufficient conditions on the nodal polynomial):
\begin{teorema}\label{22}
Let $P_n(z)$ be a polynomial of degree $n$ given by (\ref{paranewbis}) (up to a multiplicative factor). Then,
\begin{enumerate}
\item $P_n(z)$ has exactly $n$ distinct zeros $z_1,\ldots,z_n$ on $\TT$.
\item There exist positive real numbers $\lambda_1,\ldots,\lambda_n$ such that
\begin{equation}\label{nueva1}
I_n(f)=\sum_{j=1}^{n} \lambda_j f(z_j)=I_{\omega}(f)
\;\;\;,\;\;\forall f \in  \Lambda_{-(n-1),n-1}.
\end{equation}
\end{enumerate}
\begin{flushright}
$\Box$
\end{flushright}
\end{teorema}

The quadrature formula $I_n(f)$ given by (\ref{nueva1}), and earlier introduced in \cite{Jo} was called an  \lq\lq n-point Szeg\"o quadrature formula", representing the analogue on the unit circle of the Gaussian formulas for intervals of the real axis. However, it must be remarked two big differences in this respect: the nodes are not the zeros of the $n$-th orthogonal polynomial with respect to $\mu$ and the $n$-point Szeg\"o formula is exact in $\Lambda_{-(n-1),n-1}$, whose dimension is $2n-1$ instead of $2n$. Observe that since the nodes are the zeros of the $n$-th para-orthogonal polynomial characterized by (\ref{paranewbis}), an one-parameter family of quadrature formulas exact in $\Lambda_{-(n-1),n-1}$ arises.

On the other hand, starting from a generating sequence $\{ p(n) \}_{n=0}^{\infty}$ and because of the fact that $ {\cal L}_{n-1}$ is a Chebyshev system on $\TT$ of dimension $n$ (since $0 \not\in \TT$), for $n$ distinct nodes $z_1,\ldots,z_n$ on $\TT$, parameters $\lambda_1,\ldots,\lambda_n$ can be uniquely determined so that,
$I_n(L)= I_{\mu}(L)$ for all $L \in {\cal L}_{n-1}$. In order to recover Szeg\H{o} formulas in the natural framework of the orthogonal Laurent polynomials on the unit circle, and inspired by the ordinary polynomial situation, we will deal with subspaces of $\Lambda$ of the form ${\cal L}_{n} {\cal L}_{r*} = \Lambda_{-[p(n)+q(r)],[q(n)+p(r)]}$ with dimension $n+r+1$ (observe that ${\cal L}_{n-1} \subset  {\cal L}_{n} {\cal L}_{r*}$). Hence, from above, $0 \leq r \leq n-1$ and we might analyze how large $r=r(n)$ can be taken. The first step should be to consider $r=n-1$, but a negative answer is proved in \cite{RO}:

\begin{teorema}\label{import}
There cannot exist an $n$-point quadrature formula like (\ref{quad})
with nodes on $\TT$ which is exact in ${\cal L}_n {\cal L}_{(n-1)*}$
for any given arbitrary generating sequence $\{ p(n) \}_{n\geq 0}$.
\begin{flushright}
$\Box$
\end{flushright}
\end{teorema}
The second step is to consider $r=n-2$. For this purpose, we set $\lambda(n)=p(n)-p(n-2) \in \{ 0,1,2 \}$. The results obtained in \cite{RO} are summarized in:
\begin{teorema}\label{cerosquasi}
Let $\{\chi_n(z) \}_{n=0}^{\infty}$ be the sequence of orthonormal Laurent polynomials with respect to the measure $\mu$ and the ordering induced by the generating sequence $\{p(n) \}_{n=0}^{\infty}$. Suppose that $\lambda(n)=p(n)-p(n-2)=1$ and consider
\begin{equation}\label{znew1}
R_n(z,u)=C_n \left[ \eta_n \chi_n(z) + \tau_n \chi_{n-1}(z) \right]
\end{equation}
being $C_n \neq 0$ and
\begin{equation}\label{znew2}
\tau_n=\left\{ \begin{array}{ccl}
\overline{u-\delta_n} &if &s(n)=1 \\
u - \delta_n &if &s(n)=0
\end{array} \right.
\end{equation}
with $u \in \TT$ and $\{\delta_n \}_{n=0}^{\infty}$ the sequence of
Schur parameters associated with $\mu$. Then,
\begin{enumerate}
\item $R_n(z,u)$ has exactly $n$ distinct zeros on $\TT$.
\item If $z_1,\ldots,z_n$ are the zeros of $R_n(z,u)$, then there exist positive numbers $\lambda_1,\ldots,\lambda_n$ such that
\begin{equation}\label{quasiquad}
I_n(f)=\sum_{j=1}^{n} \lambda_j f(z_j) = I_{\mu}(f) \;\;\;,\;\;\forall f \in {\cal L}_{n}{\cal L}_{(n-2)*}.
\end{equation}
\item There cannot exist an $n$-point quadrature formula with nodes on $\TT$ to be exact in ${\cal L}_{n}{\cal L}_{(n-2)*}$ if $\lambda(n)=0$ or $\lambda(n)=1$.
\end{enumerate}
\begin{flushright}
$\Box$
\end{flushright}
\end{teorema}

Thus, under the assumption that $\lambda(n)=p(n)-p(n-2)=1$, we see that ${\cal L}_n {\cal L}_{(n-2)*}=\Lambda_{-(n-1),n-1}$. Therefore, the quadrature rule given by (\ref{quasiquad}) coincides with an $n$-point Szeg\H{o} quadrature formula for $\mu(\theta)$ and, taking into account that the solutions of the finite difference equation $\lambda(n)=p(n)-p(n-2)=1$ for $n \geq 2$ are given by
$$p(n)= \left\{ \begin{array}{ccl}
E\left[ \frac{n}{2} \right] &if &p(0)=p(1)=0 \\
E\left[ \frac{n+1}{2} \right] &if &p(0)=0\;,\;\; p(1)=1 \\
\end{array} \right. \;,$$
we see as the natural balanced orderings earlier introduced by Thron in \cite{Th} are again recovered. Furthermore, they are the only ones which produce quadrature formulas with nodes on $\TT$ with a maximal domain of validity. On the other hand, as we have seen in the Section 4, these orderings correspond with the narrowest matricial representation of a sequence of orthonormal Laurent polynomials.

In order to complete the construction of such quadrature formulas we give expressions for the weights also proved in \cite{RO}:
\begin{teorema}\label{pesosmejor}
Let $\{\chi_n(z) \}_{n=0}^{\infty}$ be the sequence of orthonormal Laurent polynomials with respect to the measure $\mu$ and the ordering induced by a generating sequence $\{p(n) \}_{n=0}^{\infty}$. Then, the weights $\{\lambda_j \}_{j=1}^{n}$ for the quadrature formula (\ref{quasiquad}) are given for $j=1,\ldots,n$ by, either
%\begin{enumerate}
%\item
\begin{equation}\label{newweight}
\lambda_j=\frac{1}{\sum_{k=0}^{n-1}\left| \chi_k(z_j) \right|^2}
\end{equation}
or
%\item
\begin{equation}\label{newweightbis}
\lambda_j = \frac{(-1)^{s(n)}}{2\Re \left[ z_j \chi'_n(z_j)\overline{\chi_n(z_j)}\right] + (p(n)-q(n))\left| \chi_n(z_j) \right|^2} \;,
\end{equation}
%\end{enumerate}
where the nodes $\{z_j\}_{j=1}^{n}$ are the zeros of $R_n(z,u)$ given by (\ref{znew1}), or equivalently the zeros of $P_n(z)$ in (\ref{paranewbis}).
\begin{flushright}
$\Box$
\end{flushright}
\end{teorema}

Once quadrature formulas on $\TT$ with a maximal domain of validity have been constructed (nodes from (\ref{paranewbis}) or (\ref{znew1}) and weights by (\ref{newweight}) or (\ref{newweightbis})) it seems that the following should be done is their effective computation. In this respect, it should be noticed that, untill the present moment and as far as we know, numerical experiments with Szeg\H{o} quadrature have mostly involved either measures whose sequences of Szeg\H{o} polynomials are explicitly known or these polynomials have been computed by Levinson's algorithm (see e.g. \cite{Le}, \cite{LD} and \cite{JC}). In this respect, the zeros of (\ref{paranewbis}) or (\ref{znew1}) can be found by using any standard root finding method available in the literature (as for a specific procedure concerning rational modifications of the Lebesgue measure see also \cite{Tru}). Hence, in the rest of the section we shall revise some strategies to effectively compute the nodes $\{ z_j \}_{j=1}^{n}$ and weights $\{ \lambda_j \}_{j=1}^{n}$ for an $n$-point Szeg\H{o} quadrature formula:
$$I_n(f)=\sum_{j=1}^{n} \lambda_j f(z_j) \;\;\;,\;\;|z_j|=1 \;\;\;,\;\;j=1,\ldots,n.$$
Indeed, for the nodal polynomial $P_n(z)=\prod_{j=1}^{n} \left( z-z_j \right)$ given by (\ref{paranewbis}) we can write from (\ref{le}) (see \cite{CM}):
$$P_n(z)=P_n(z,u)=z\rho_{n-1}(z) + u\rho_{n-1}^*(z) \;\;,\;|u|=1.$$
Now, it can be easily checked that $\{z\varphi_0(z),\ldots,z\varphi_{n-2}(z),-u\varphi_{n-1}^*(z) \}$ is an orthonormal basis of $\PP_{n-1}$ which must be related to $\{\varphi_0(z),\ldots,\varphi_{n-1}(z) \}$ by an unitary matrix $U_n$. Setting for $n\geq 0$, $e_n=\langle \rho_n(z),\rho_n(z) \rangle_{\mu}$ (recall that $e_n=\prod_{k=1}^{n} \eta_k^2$ for $n \geq 1$ and that $\varphi_n(z)=\frac{\rho_n(z)}{\sqrt{e_n}}$) it is straight forward to check that $U_n$ is the following irreducible Hessenberg matrix:
\begin{equation}\label{un}
U_n= \left( \begin{array}{ccccc}
d_{1,1} &d_{1,2} &0 &\cdots &0 \\
d_{2,1} &d_{2,2} &d_{2,3} &\cdots &0 \\
\vdots &\vdots &\vdots &\ddots &\vdots \\
d_{n,1} &d_{n,2} &d_{n,3} &\cdots &d_{n,n}
\end{array} \right)
\end{equation}
where
\begin{equation}\label{undata}
d_{i,j}= \left\{ \begin{array}{crl}
-\overline{\delta_{j-1}}u\sqrt{\frac{e_{n-1}}{e_{j-1}}} &if &i=n \\
-\overline{\delta_{j-1}}\delta_{i}\sqrt{\frac{e_{i-1}}{e_{j-1}}} &if &i \leq n-1 \;\;,\;j \leq i \\
\eta_i &if &i \leq n-1 \;\;,\;j=i+1
\end{array} \right. .
\end{equation}
Hence, we can write
$$\left( \begin{array}{c}
z\varphi_0(z) \\
\vdots \\
z\varphi_{n-2}(z) \\
-u\varphi_{n-1}^*(z)
\end{array} \right) = U_n \left( \begin{array}{c}
\varphi_0(z) \\
\vdots \\
\varphi_{n-2}(z) \\
\varphi_{n-1}(z)
\end{array} \right) .$$
Observe that, if we introduce the notation ${\cal H}_n(\delta_0,\ldots,\delta_{n})$ to indicate the principal submatrix of ${\cal H}(\delta)$ (given by (\ref{hess})) of order $n$ with entries $\{ \delta_k \}_{k=0}^{n}$, then $U_n={\cal H}_n(\delta_0,\ldots,\delta_{n-1},u)$, that is, the principal submatrix of ${\cal H}(\delta)$ of order $n$ where the Schur parameter $\delta_n$ is replaced by $u \in \TT$. Now,
$$\left( \begin{array}{c}
z\varphi_0(z) \\
\vdots \\
z\varphi_{n-2}(z) \\
-u\varphi_{n-1}^*(z)
\end{array} \right) = \left( \begin{array}{c}
z\varphi_0(z) \\
\vdots \\
z\varphi_{n-2}(z) \\
z\varphi_{n-1}(z)-z\varphi_{n-1}(z)-u\varphi_{n-1}^*(z)
\end{array} \right) = z \left( \begin{array}{c}
\varphi_0(z) \\
\vdots \\
\varphi_{n-2}(z) \\
\varphi_{n-1}(z)
\end{array} \right) - \left( \begin{array}{c}
0 \\
\vdots \\
0 \\
\frac{P_n(z,u)}{\sqrt{e_{n-1}}}
\end{array} \right).$$
Thus, we have the identity
\begin{equation}\label{zvn}
zV_n(z)=U_nV_n(z) + b_n(z)
\end{equation}
where
$$V_n(z)=\left( \varphi_0(z), \varphi_1(z), \cdots,\varphi_{n-1}(z) \right)^{T} \;\;;\;\;\;b_n(z)=\left( 0,\cdots,0,\frac{P_n(z,u)}{\sqrt{e_{n-1}}} \right)^T$$
with $|u|=1$. From (\ref{zvn}) one sees that any zero $\xi$ of
$P_n(z,u)$ is an eigenvalue of $U_n$ with associated eigenvector
$V_n(\xi)$. So, let $z_j$ be a zero of $P_n(z,u)$ for $j=1,\ldots,n$
and consider the corresponding normalized eigenvector of $U_n$:
$$W_n(z_j)=\frac{V_n(z_j)}{\left[ \sum_{k=0}^{n-1} \left| \varphi_k(z_j) \right|^2 \right]^{1/2}}.$$
Now, taking into account that for any $z \in \TT$, $\left| \varphi_k(z) \right|^2=\left| \chi_k(z) \right|^2$, from (\ref{newweight}) it follows that
\begin{equation}\label{wn}
W_n(z_j)=\lambda_j^{1/2}V_n(z_j).
\end{equation}
If we write $W_n(z_j)=\left( q_{0,j},\ldots,q_{n-1,j} \right)^T$ and select out the first components of both sides of (\ref{wn}), we obtain $q_{0,j}=\lambda_j^{1/2}\varphi_0(z_j)$. But, since we are dealing with a probability measure ($\int d\mu(\theta)=1$) then $\varphi_0(z) \equiv 1$ and hence,
$$\lambda_j=q_{0,j}^2 \;\;\;,\;\;j=1,\ldots,n.$$
In short, the following theorem has been proved

\begin{teorema}\label{eigen1}
Let $I_n(f)$ be the $n$-th Szeg\H{o} quadrature formula (\ref{quasiquad}). Then,
\begin{enumerate}
\item The nodes $\{z_j \}_{j=1}^{n}$ are the eigenvalues of $U_n={\cal H}_n(\delta_0,\ldots,\delta_{n-1},u)$ given by (\ref{un})-(\ref{undata}), for all $u \in \TT$.
\item The weights $\{\lambda_j \}_{j=1}^{n}$ are given by the first component of the normalized eigenvectors.
\end{enumerate}
\end{teorema}

Finally, we show an alternative approach to the computation of a Szeg\H{o} quadrature formula (\ref{quasiquad}) by using truncations of the five-diagonal matrix ${\cal C}(\delta)$. In the next result we will use the matrix ${\cal C}_n(\delta_1,\ldots,\delta_{n-1},u)$, that is, the $n$-th principal submatrix of ${\cal C}(\delta)$ of order $n$ where the Schur parameter $\delta_n$ is replaced by $u \in \TT$. The first part has been already deduced in \cite{CMV4} by using operator theory techniques meanwhile our proof here presented is based on the recurrence relations satisfied by the family of orthonormal Laurent polynomials. Without loss of generality, we can fix the ordering induced by $p(n)=E\left[ \frac{n}{2} \right]$ (recall that the matrix representation associated with the ordering induced by $p(n)=E\left[ \frac{n+1}{2} \right]$ is ${\cal C}(\delta)^T$).

\begin{teorema}\label{eigen2}
Let $I_n(f)$ be the $n$-th Szeg\H{o} quadrature formula (\ref{quasiquad}). Then,
\begin{enumerate}
\item The nodes $\{z_j \}_{j=1}^{n}$ are the eigenvalues of ${\cal C}_n(\delta_1,\ldots,\delta_{n-1},u)$, for all $u \in \TT$.
\item The weights $\{\lambda_j \}_{j=1}^{n}$ are given by the first component of the normalized eigenvectors.
\end{enumerate}
\end{teorema}

{\em Proof}.- As we have seen, the nodes are the zeros of $R_n(z)$ given by (\ref{znew1})-(\ref{znew2}). Setting $X_n(z)= \left( \chi_0(z),\ldots,\chi_{n-1}(z)\right)^T$ it follows from (\ref{five2}) that
$$zX_n(z)={\cal C}_n(\delta_1,\ldots,\delta_{n-1},u) X_n(z) + T_n(z) \;,$$
where
$$T_n(z)= \left\{ \begin{array}{ccl}
(0,\ldots,0,A(z))^T &if &n=2k \\
(0,\ldots,0,B(z),C(z))^T &if &n=2k+1
\end{array} \right.$$
being
$$\begin{array}{rcl}
A(z) &=
&\eta_{2k-1}(u-\delta_{2k})\chi_{2k-2}(z)+\overline{\delta_{2k-1}}(u-\delta_{2k})\chi_{2k-1}(z)-
\\
&&\eta_{2k}\delta_{2k+1}\chi_{2k}(z)+\eta_{2k}\eta_{2k+1}\chi_{2k+1}(z)
\end{array}$$
$$B(z)=\eta_{2k}(u-\delta_{2k+1})\chi_{2k}(z) + \eta_{2k}\eta_{2k+1}\chi_{2k+1}(z)$$
$$C(z)=\overline{\delta_{2k}}(u-\delta_{2k+1})\chi_{2k}(z) + \overline{\delta_{2k}}\eta_{2k+1}\chi_{2k+1}(z).$$
Now, on the one hand, in the odd case ($n=2k+1$) it follows from (\ref{znew1})-(\ref{znew2}) that
$$T_n(z)=\left( 0,\ldots,0, \eta_{n-1}R_n(z),\overline{\delta_{n-1}}R_n(z) \right)^T$$
implying that $T_n(z) \equiv (0,\ldots,0)^T$ if and only if $R_n(z) =0$. On the other hand, from (\ref{eq4}) and since $u \in \TT$ it follows in the even case ($n=2k$) that
$$\begin{array}{ccl}
A(z) &= &(u-\delta_{2k})\left[ z\eta_{2k}\chi_{2k}(z)-\left( \overline{\delta_{2k}}z + \overline{\delta_{2k-1}} \right) \chi_{2k-1}(z) \right]\\
&& + \overline{\delta_{2k-1}}(u-\delta_{2k})\chi_{2k-1}(z)-\eta_{2k}\delta_{2k+1}\chi_{2k}(z) + \\
&& \eta_{2k}\left[ \left( \delta_{2k+1} + \delta_{2k}z
\right)\chi_{2k}(z) + \eta_{2k}z\chi_{2k-1}(z) \right] \\
&= &uzR_n(z)
\end{array}$$
and hence it holds again that $T_n(z) \equiv (0,\ldots,0)^T$ if and only if $R_n(z)=0$. The proof of the second part follows from the same arguments as in Theorem \ref{eigen1}.
\begin{flushright}
$\Box$
\end{flushright}

%In \cite{CMV4} it is proved that for all $u \in \TT$, the eigenvalues of ${\cal C}_n(\delta_1,\ldots,\delta_{n-1},u)$ are simple and coincides with the zeros of $T_n(z,u)=z\varphi_{n-1}(z)+u\varphi_{n-1}^*(z)$. Fixed $p(n)=E\left[ \frac{n+1}{2} \right]$ it follows from Theorem \ref{cerosquasi} and (\ref{lenorma}) that the nodes of the quadrature formula (\ref{quasiquad}) are, for all $\alpha_n \in \TT$, the zeros of
%$$R_n(z)=\left\{ \begin{array}{lcl}
%z^{-k}T_{2k}(z,\alpha_{2k}) &if &n=2k \\
%\\
%\overline{\alpha_{2k+1}}z^{-(k+1)}T_{2k+1}(z,\alpha_{2k+1}) &if &n=2k+1
%\end{array} \right. ,$$
%leading with an alternative approach to compute the nodes of quadrature formulas as eigenvalues of truncations of ${\cal C}(\delta)$.

\section{Numerical examples}

\setcounter{equation}{0}

In order to numerically illustrate the results given in the previous section we expose the computation of nodes and weights of the Szeg\H{o} quadrature formulas considering the absolutely continuous measure defined on $[-\pi,\pi]$ by $d\mu(\theta)=\omega(\theta)d\theta$ with
\begin{equation}\label{rogerszegoweight}
\omega(\theta) = \frac{2\pi}{\sqrt{2\pi \log \left( \frac{1}{q} \right)}} \sum_{j=-\infty}^{\infty} \exp \left({\frac{-\left(\theta - 2\pi j \right)^2}{2\log \left( \frac{1}{q} \right)}} \right) \;\;\;,\;\;q \in (0,1) \;\;.
\end{equation}
The corresponding monic orthogonal polynomials are the so-called Rogers-Szeg\"o $q$-polynomials. Throughout this section we also fix the ordering induced by the generating sequence $p(n)=E \left[ \frac{n+1}{2} \right]$. An explicit expression for such polynomials is given in \cite[Ch. 1]{BS} and so from Proposition \ref{conex2} the following explicit expression for the corresponding monic orthogonal Laurent polynomials is deduced:
\begin{equation}\label{rogsez}
\phi_n(z)= \left\{ \begin{array}{lcl}
\sum_{j=-k}^{k} (-1)^{j+k} \left[
\begin{array}{c}
2k \\
j+k
\end{array}
\right]_q q^{\frac{k-j}{2}} z^{j} &if &n=2k \\
\sum_{j=-(k+1)}^{k} (-1)^{j+k+1} \left[
\begin{array}{c}
2k+1 \\
k-j
\end{array}
\right]_q q^{\frac{j+k+1}{2}} z^{j} &if &n=2k+1
\end{array} \right.
\end{equation}
where, as usual, the $q$-binomial coefficients $\left[
\begin{array}{c}
n \\
j
\end{array}
\right]_q$ are defined by,
\begin{equation}
\begin{array}{ccl}
(n)_q &= &(1-q)(1-q^2)\cdots (1-q^n) \;\;\;,\;\;(0)_q \equiv 1 \\
\\
\left[
\begin{array}{c}
n \\
j
\end{array}
\right]_q &= &\frac{(n)_q}{(j)_q (n-j)_q} = \frac{(1-q^n) \cdots (1-q^{n-j+1})}{(1-q) \cdots (1-q^j)}.
\end{array}
\end{equation}
Now, writting $\phi_{2k}(z)=\sum_{j=-k}^{k} a_jz^j$ and $\phi_{2k+1}(z)=\sum_{j=-(k+1)}^{k} b_jz^j$, then the coefficients $\{a_j \}_{j=-k}^{k}$ and $\{b_j \}_{j=-(k+1)}^{k}$ can be recursively computed by
$$a_k=1 \;\;\;,\;\;\;a_j=-a_{j+1}\sqrt{q}\frac{1-q^{k+j+1}}{1-q^{k-j}} \;\;,\;-k \leq j \leq k-1$$
$$b_{-(k+1)}=1 \;\;\;,\;\;\;b_{j+1}=-b_j\sqrt{q}\frac{1-q^{k-j}}{1-q^{k+j+2}} \;\;,\;-(k+1) \leq j \leq k-1.$$

The following two figures show the zeros of $\phi_{10}(z)$ and $\phi_{11}(z)$, taking $q=0.1, 0.25, 0.5, 0.75$ and $0.9$ (as usual, the circles below represent the unit circle). The data were computed by using the Jenkins and Traub root-finding method (see \cite{JTR}), appropriate in this case since the polynomials has real coefficients, and they are similar to those given in \cite[Ch. 8]{BS} where the distribution of zeros of Szeg\"o polynomials for this measure is also considered:

$$Figure \;1 \;\;\left( zeros \;of \;\phi_{10}(z) \right)$$
$$\begin{array}{ccc}
q=0.1 & q=0.25 &q=0.5 \\
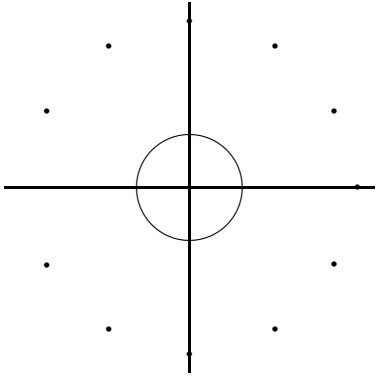
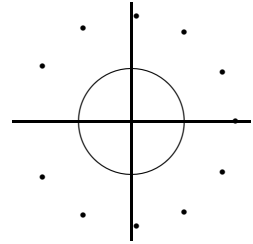
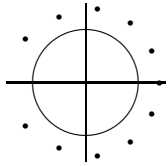
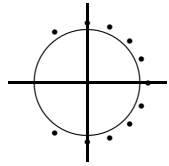
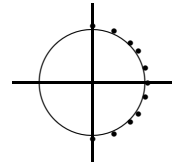
\begin{picture}(120,70)
\put(60,35){\circle{60}}
\put(30,35){\line(1,0){60}}
\put(60,5){\line(0,1){60}}
\put(55,38.5){\circle*{2}}
\put(55,31.5){\circle*{2}}
\put(57.5,40.25){\circle*{2}}
\put(57.5,29.75){\circle*{2}}
\put(60,41){\circle*{2}}
\put(60,29){\circle*{2}}
\put(64.25,39.5){\circle*{2}}
\put(64.25,30.5){\circle*{2}}
\put(66.5,35.25){\circle*{2}}
\put(66.5,34.75){\circle*{2}}
\end{picture} &
\begin{picture}(120,70)
\put(60,35){\circle{60}}
\put(30,35){\line(1,0){60}}
\put(60,5){\line(0,1){60}}
\put(51,40){\circle*{2}}
\put(51,30){\circle*{2}}
\put(55.75,43.75){\circle*{2}}
\put(55.75,26.25){\circle*{2}}
\put(61.5,45){\circle*{2}}
\put(61.5,25){\circle*{2}}
\put(66.5,42){\circle*{2}}
\put(66.5,28){\circle*{2}}
\put(69,37.5){\circle*{2}}
\put(69,32.5){\circle*{2}}
\end{picture} &
\begin{picture}(120,70)
\put(60,35){\circle{60}}
\put(30,35){\line(1,0){60}}
\put(60,5){\line(0,1){60}}
\put(48.5,44){\circle*{2}}
\put(48.5,26){\circle*{2}}
\put(56.5,48.5){\circle*{2}}
\put(56.5,21.5){\circle*{2}}
\put(63.75,48.25){\circle*{2}}
\put(63.75,21.75){\circle*{2}}
\put(70.5,44){\circle*{2}}
\put(70.5,26){\circle*{2}}
\put(73.25,38.5){\circle*{2}}
\put(73.25,31.5){\circle*{2}}
\end{picture}
\end{array}$$
$$\begin{array}{cc}
q=0.75 & q=0.9 \\
\begin{picture}(120,70)
\put(60,35){\circle{60}}
\put(30,35){\line(1,0){60}}
\put(60,5){\line(0,1){60}}
\put(52,49.5){\circle*{2}}
\put(52,20.5){\circle*{2}}
\put(60,52.25){\circle*{2}}
\put(60,17.75){\circle*{2}}
\put(68,50){\circle*{2}}
\put(68,20){\circle*{2}}
\put(74,45){\circle*{2}}
\put(74,25){\circle*{2}}
\put(77,39){\circle*{2}}
\put(77,31){\circle*{2}}
\end{picture} &
\begin{picture}(120,70)
\put(60,35){\circle{60}}
\put(30,35){\line(1,0){60}}
\put(60,5){\line(0,1){60}}
\put(61.5,53.5){\circle*{2}}
\put(61.5,16.5){\circle*{2}}
\put(68.5,51.5){\circle*{2}}
\put(68.5,18.5){\circle*{2}}
\put(73.5,47.5){\circle*{2}}
\put(73.5,22.5){\circle*{2}}
\put(77.15,41.6){\circle*{2}}
\put(77.15,28.4){\circle*{2}}
\put(78.5,37){\circle*{2}}
\put(78.5,33){\circle*{2}}
\end{picture}
\end{array}$$

%%%%%%%%%%%%%%%%%%%%%%%%%%%%%%%%%%%%%%%%%%%%%%%%%%%%%%%%%%%%%%%%%%%%%%%%%%%%%%%%%%ELIMINATE!!
%\newpage
%%%%%%%%%%%%%%%%%%%%%%%%%%%%%%%%%%%%%%%%%%%%%%%%%%%%%%%%%%%%%%%%%%%%%%%%%%%%%%%%%%%%%%%%%%%%%

$$Figure \;2 \;\;\left( zeros \;of \;\phi_{11}(z) \right)$$
$$\begin{array}{lr}
q=0.1 & q=0.25 \\
\\
\\
\\
\\
\begin{picture}(120,70)
\put(15,35){\circle{60}}
\put(-55,35){\line(1,0){140}}
\put(15,-35){\line(0,1){140}}
\put(78.5,35){\circle*{2}}
\put(69.75,64){\circle*{2}}
\put(69.75,6){\circle*{2}}
\put(47.5,88.5){\circle*{2}}
\put(47.5,-18.5){\circle*{2}}
\put(15.1,98){\circle*{2}}
\put(15.1,-28){\circle*{2}}
\put(-15.5,88.6){\circle*{2}}
\put(-15.5,-18.6){\circle*{2}}
\put(-39,64.25){\circle*{2}}
\put(-39,5.75){\circle*{2}}
\end{picture} &
\begin{picture}(120,70)
\put(105,35){\circle{60}}
\put(60,35){\line(1,0){90}}
\put(105,-10){\line(0,1){90}}
\put(144.5,35){\circle*{2}}
\put(139.5,54){\circle*{2}}
\put(139.5,16){\circle*{2}}
\put(125,69){\circle*{2}}
\put(125,1){\circle*{2}}
\put(107,74.75){\circle*{2}}
\put(107,-4.75){\circle*{2}}
\put(86.75,70.5){\circle*{2}}
\put(86.75,-0.5){\circle*{2}}
\put(71.5,56){\circle*{2}}
\put(71.5,14){\circle*{2}}
\end{picture}
\end{array}$$
\vspace{2cm}
$$\begin{array}{ccc}
q=0.5 & q=0.75 &q=0.9 \\
\\
\begin{picture}(120,70)
\put(60,35){\circle{60}}
\put(30,35){\line(1,0){60}}
\put(60,5){\line(0,1){60}}
\put(88,35){\circle*{2}}
\put(85,47){\circle*{2}}
\put(85,23){\circle*{2}}
\put(77,57.5){\circle*{2}}
\put(77,12.5){\circle*{2}}
\put(50,59.75){\circle*{2}}
\put(50,10.25){\circle*{2}}
\put(64.5,62.75){\circle*{2}}
\put(64.5,7.25){\circle*{2}}
\put(37.25,51.5){\circle*{2}}
\put(37.25,18.5){\circle*{2}}
\end{picture} &
\begin{picture}(120,70)
\put(60,35){\circle{60}}
\put(30,35){\line(1,0){60}}
\put(60,5){\line(0,1){60}}
\put(83,35){\circle*{2}}
\put(80.5,44){\circle*{2}}
\put(80.5,26){\circle*{2}}
\put(76,50.75){\circle*{2}}
\put(76,19.25){\circle*{2}}
\put(60,57.5){\circle*{2}}
\put(60,12.5){\circle*{2}}
\put(68.5,55.9){\circle*{2}}
\put(68.5,14.1){\circle*{2}}
\put(47.75,54){\circle*{2}}
\put(47.75,16){\circle*{2}}
\end{picture} &
\begin{picture}(120,70)
\put(60,35){\circle{60}}
\put(30,35){\line(1,0){60}}
\put(60,5){\line(0,1){60}}
\put(81,35){\circle*{2}}
\put(80.25,40.5){\circle*{2}}
\put(80.25,29.5){\circle*{2}}
\put(77.5,46.75){\circle*{2}}
\put(77.5,23.25){\circle*{2}}
\put(74.5,50){\circle*{2}}
\put(74.5,20){\circle*{2}}
\put(68.25,54.4){\circle*{2}}
\put(68.25,15.6){\circle*{2}}
\put(60.25,56.4){\circle*{2}}
\put(60.25,13.6){\circle*{2}}
\end{picture}
\end{array}$$
From Figures 1 and 2 one can observe that the zeros of $\phi_{10}(z)$ are located on the circle $\left\{z : |z|=q^{1/2} \right\}$ whereas the zeros of $\phi_{11}(z)$ are located on the circle $\left\{z : |z|=q^{-1/2} \right\}$, in accordance with Proposition \ref{conex2} and the Mazel-Geronimo-Hayes theorem (see \cite{MZ}).

Now, we compute the nodes and weights of the quadrature formula (\ref{quasiquad}) also for $n=10$ and $q=0.1, 0.25, 0.5, 0.75$ and $0.9$. The nodes are the zeros of $R_{10}(z,u)$ given by (\ref{znew1}) (we take $u=1$). The corresponding sequence of Schur parameters are given by $\delta_n=(-1)^nq^{\frac{n}{2}}$ for all $n=1,2,\ldots$ (see \cite{BS}). As for the weights, we can make use of (\ref{newweightbis}) taking into account that
$$\chi_n(z)=\frac{\phi_n(z)}{\sqrt{(1-q)\cdots (1-q^n)}}$$
with $\phi_n(z)$ explicitly given by (\ref{rogsez}). The results are displayed in Tables 1-5.

%%%%%%%%%%%%%%%%%%%%%%%%%%%%%%%%%%%%%%%%%%%%%%%%%%%%%%%%%%%%%%%%%%%%%%%%%%%%%%%%%%ELIMINATE!!
\newpage
%%%%%%%%%%%%%%%%%%%%%%%%%%%%%%%%%%%%%%%%%%%%%%%%%%%%%%%%%%%%%%%%%%%%%%%%%%%%%%%%%%%%%%%%%%%%%

%\begin{small}
\begin{center}
$Table \;1 \;\;\;\;\;\;\;\;\;\;\;\;\;\;\;\;\;\;\;\;\;\;\;\;\;\;\;\;\;\;\;\;\;\;\;\;\;\;\;\;\;\;\;\;\;\;\;\;\;\;\;\;\;\;\;\;Table \;2$
\end{center}
$$\begin{array}{rl}
\hfill
\hbox{
\small\tabcolsep 5pt
\begin{tabular}{|c|c|}
\hline & \\
$q=0.1$ & $n=10$ \\
\hline & \\
Nodes & Weights \\
\hline & \\
$\begin{array}{r}
-0.940400 \pm 0.34007i \\
-0.531157 \pm 0.847273i \\
0.0668824 \pm 0.997761i \\
0.624424 \pm 0.781086i \\
0.955949 \pm 0.293533i
\end{array}$ &
$\begin{array}{c}
0.0459602 \\
0.0669775 \\
0.100057 \\
0.133157 \\
0.153848
\end{array}$ \\
& \\
\hline
\end{tabular}}
\hfill\bigskip\par\noindent
&
\hfill
\hbox{
\small\tabcolsep 5pt
\begin{tabular}{|c|c|}
\hline & \\
$q=0.25$ & $n=10$ \\
\hline & \\
Nodes & Weights \\
\hline & \\
$\begin{array}{r}
-0.922051 \pm 0.387069i \\
-0.473103 \pm 0.881007i \\
0.119954 \pm 0.992779i \\
0.650270 \pm 0.759703i \\
0.959239 \pm 0.282596i
\end{array}$ &
$\begin{array}{c}
0.0195773 \\
0.0466391 \\
0.0947585 \\
0.150362 \\
0.188665
\end{array}$ \\
& \\
\hline
\end{tabular}}
\hfill\bigskip\par\noindent
\end{array}$$

\begin{center}
$Table \;3 \;\;\;\;\;\;\;\;\;\;\;\;\;\;\;\;\;\;\;\;\;\;\;\;\;\;\;\;\;\;\;\;\;\;\;\;\;\;\;\;\;\;\;\;\;\;\;\;\;\;\;\;\;\;\;\;Table \;4$
\end{center}
$$\begin{array}{rl}
\hfill
\hbox{
\small\tabcolsep 5pt
\begin{tabular}{|c|c|}
\hline & \\
$q=0.5$ & $n=10$ \\
\hline & \\
Nodes & Weights \\
\hline & \\
$\begin{array}{r}
-0.842988 \pm 0.537932i \\
-0.333209 \pm 0.942853i \\
0.234605 \pm 0.972091i \\
0.703537 \pm 0.710659i \\
0.965879 \pm 0.258994i
\end{array}$ &
$\begin{array}{c}
0.00312009 \\
0.0207928 \\
0.0737936 \\
0.163017 \\
0.239274
\end{array}$ \\
& \\
\hline
\end{tabular}}
\hfill\bigskip\par\noindent
&
\hfill
\hbox{
\small\tabcolsep 5pt
\begin{tabular}{|c|c|}
\hline & \\
$q=0.75$ & $n=10$ \\
\hline & \\
Nodes & Weights \\
\hline & \\
$\begin{array}{r}
-0.517559 \pm 0.855648i \\
0.00961854 \pm 0.999954i \\
0.467501 \pm 0.883993i \\
0.801825 \pm 0.597559i \\
0.977622 \pm 0.210369i
\end{array}$ &
$\begin{array}{c}
0.000196919 \\
0.00541542 \\
0.0439839 \\
0.158275 \\
0.292128
\end{array}$ \\
& \\
\hline
\end{tabular}}
\hfill\bigskip\par\noindent
\end{array}$$

%%%%%%%%%%%%%%%%%%%%%%%%%%%%%%%%%%%%%%%%%%%%%%%%%%%%%%%%%%%%%%%%%%%%%%%%%%%%%%%%%%ELIMINATE!!
\newpage
%%%%%%%%%%%%%%%%%%%%%%%%%%%%%%%%%%%%%%%%%%%%%%%%%%%%%%%%%%%%%%%%%%%%%%%%%%%%%%%%%%%%%%%%%%%%%

\begin{center}
$Table \;5$
\end{center}
%$$\begin{array}{rl}
\hfill
\hbox{
\small\tabcolsep 5pt
\begin{tabular}{|c|c|}
\hline & \\
$q=0.9$ & $n=10$ \\
\hline & \\
Nodes & Weights \\
\hline & \\
$\begin{array}{r}
0.112467 \pm 0.993655i \\
0.475746 \pm 0.879582i \\
0.734593 \pm 0.678508i \\
0.904709 \pm 0.426031i \\
0.989421 \pm 0.145076i
\end{array}$ &
$\begin{array}{c}
0.0000221869 \\
0.000173911 \\
0.0276214 \\
0.146351 \\
0.324221
\end{array}$ \\
& \\
\hline
\end{tabular}}
\hfill\bigskip\par\noindent
%\end{small}

\begin{nota}
From the above tables one sees that the weights corresponding to a pair of complex conjugate nodes are equal. This directly follows from (\ref{newweight}) since in this case the coefficients of the $n$-th orthonormal Laurent polynomial are real.
\end{nota}

In order to check the effectiveness of Theorems \ref{eigen1} and \ref{eigen2} we also compute for $n=10$, $q=0.1,0.25,0.5,0.75,0.9$ and $u=1$ the eigenvalues and the first component of the normalized eigenvectors of ${\cal C}_{10}(\delta_1,\ldots,\delta_9,1)$ and ${\cal H}_{10}(\delta_0,\ldots,\delta_9,1)$, yielding the nodes and weights of the quadrature formula (\ref{quasiquad}). The computations were made by using a standard eigenvalue-finding method and the results were exactly the same ones as those obtained before. Moreover, the results of the computations for $n=11$ are also showed in the following tables:

%%%%%%%%%%%%%%%%%%%%%%%%%%%%%%%%%%%%%%%%%%%%%%%%%%%%%%%%%%%%%%%%%%%%%%%%%%%%%%%%%%ELIMINATE!!
\newpage
%%%%%%%%%%%%%%%%%%%%%%%%%%%%%%%%%%%%%%%%%%%%%%%%%%%%%%%%%%%%%%%%%%%%%%%%%%%%%%%%%%%%%%%%%%%%%

\begin{center}
$Table \;6 \;\;\;\;\;\;\;\;\;\;\;\;\;\;\;\;\;\;\;\;\;\;\;\;\;\;\;\;\;\;\;\;\;\;\;\;\;\;\;\;\;\;\;\;\;\;\;\;\;\;\;\;\;\;\;\;Table \;7$
\end{center}
$$\begin{array}{rl}
\hfill
\hbox{
\small\tabcolsep 5pt
\begin{tabular}{|c|c|}
\hline & \\
$q=0.1$ & $n=11$ \\
\hline & \\
Nodes & Weights \\
\hline & \\
$\begin{array}{r}
-1 \\
-0.814432 \pm 0.580259i \\
-0.355944 \pm 0.934507i \\
0.199254 \pm 0.979948i \\
0.68359 \pm 0.729866i \\
0.963208 \pm 0.268755i
\end{array}$ &
$\begin{array}{c}
0.03873768 \\
0.04723999 \\
0.06942034 \\
0.09817396 \\
0.12485883 \\
0.14093668
\end{array}$ \\
& \\
\hline
\end{tabular}}
\hfill\bigskip\par\noindent
&
\hfill
\hbox{
\small\tabcolsep 5pt
\begin{tabular}{|c|c|}
\hline & \\
$q=0.25$ & $n=11$ \\
\hline & \\
Nodes & Weights \\
\hline & \\
$\begin{array}{r}
-1 \\
-0.778821 \pm 0.627246i \\
-0.301023 \pm 0.953617i \\
0.243118 \pm 0.969997i \\
0.70384 \pm 0.710359i \\
0.965731 \pm 0.259545i
\end{array}$ &
$\begin{array}{c}
0.01403402 \\
0.00025263 \\
0.05384075 \\
0.09739723 \\
0.14346753 \\
0.17367845
\end{array}$ \\
& \\
\hline
\end{tabular}}
\hfill\bigskip\par\noindent
\end{array}$$

\begin{center}
$Table \;8 \;\;\;\;\;\;\;\;\;\;\;\;\;\;\;\;\;\;\;\;\;\;\;\;\;\;\;\;\;\;\;\;\;\;\;\;\;\;\;\;\;\;\;\;\;\;\;\;\;\;\;\;\;\;\;\;Table \;9$
\end{center}
$$\begin{array}{rl}
\hfill
\hbox{
\small\tabcolsep 5pt
\begin{tabular}{|c|c|}
\hline & \\
$q=0.5$ & $n=11$ \\
\hline & \\
Nodes & Weights \\
\hline & \\
$\begin{array}{r}
-1 \\
-0.681725 \pm 0.731608i \\
-0.178603 \pm 0.983921i \\
0.335240 \pm 0.942133i \\
0.745084 \pm 0.666970i \\
0.970795 \pm 0.239909i
\end{array}$ &
$\begin{array}{c}
0.00105402 \\
0.00606822 \\
0.02834884 \\
0.08185750 \\
0.16052782 \\
0.22291291
\end{array}$ \\
& \\
\hline
\end{tabular}}
\hfill\bigskip\par\noindent
&
\hfill
\hbox{
\small\tabcolsep 5pt
\begin{tabular}{|c|c|}
\hline & \\
$q=0.75$ & $n=11$ \\
\hline & \\
Nodes & Weights \\
\hline & \\
$\begin{array}{r}
-1 \\
-0.400355 \pm 0.91636i \\
0.0972837 \pm 0.995257i \\
0.517853 \pm 0.85547i \\
0.821232 \pm 0.570595i \\
0.979848 \pm 0.199745i
\end{array}$ &
$\begin{array}{c}
0.00000331 \\
0.00047627 \\
0.00836960 \\
0.05348059 \\
0.16811708 \\
0.29160532
\end{array}$ \\
& \\
\hline
\end{tabular}}
\hfill\bigskip\par\noindent
\end{array}$$

%%%%%%%%%%%%%%%%%%%%%%%%%%%%%%%%%%%%%%%%%%%%%%%%%%%%%%%%%%%%%%%%%%%%%%%%%%%%%%%%%%ELIMINATE!!
\newpage
%%%%%%%%%%%%%%%%%%%%%%%%%%%%%%%%%%%%%%%%%%%%%%%%%%%%%%%%%%%%%%%%%%%%%%%%%%%%%%%%%%%%%%%%%%%%%

\begin{center}
$Table \;10$
\end{center}
%$$\begin{array}{rl}
\hfill
\hbox{
\small\tabcolsep 5pt
\begin{tabular}{|c|c|}
\hline & \\
$q=0.9$ & $n=11$ \\
\hline & \\
Nodes & Weights \\
\hline & \\
$\begin{array}{r}
-1 \\
0.149722 \pm 0.988728i \\
0.498881 \pm 0.866671i \\
0.746631 \pm 0.665238i \\
0.909098 \pm 0.416581i \\
0.989911 \pm 0.141688i
\end{array}$ &
$\begin{array}{c}
0.00000000058 \\
0.00005234 \\
0.00328985 \\
0.04447931 \\
0.21730036 \\
0.46353159
\end{array}$ \\
& \\
\hline
\end{tabular}}
\hfill\bigskip\par\noindent

Here, it should be remarked that when computations for higher values of $n$ are required, special eigenvalue-finding methods should be considered because of the error propagation. In this respect, appropriate procedures for product of matrices (see \cite{W}) could be considered due to the factorization of ${\cal C}(\delta)$ as product of two tri-diagonal matrices. Moreover, in \cite{Ca} it is remarked that:
\begin{enumerate}
\item The computational cost by using techniques for banded matrices in comparision with techniques for Hessenberg matrices is reduced.
\item In the ${\cal C}(\delta)$ matrix, the Schur parameters appear in only finitely many elements and hence, every modification of a finite number of Schur parameters induces a finitely dimensional perturbation, something that is not true in the Hessenberg matrix.
\end{enumerate}
Throughout this section the computations were made using MATHEMATICA software.

\vspace{0.5cm}

{\bf Aknowledgements}

The authors are very grateful to Professor Olav Nj\aa stad from the University of Trondheim (Norway) for his remarks and useful sugestions.

\end{document}